\documentclass[12pt]{article}

\usepackage{amsmath}
\usepackage{amssymb}
\usepackage{theorem}
\theorembodyfont{\rmfamily}
\usepackage[mathscr]{eucal}
\newtheorem{theorem}{Theorem}[section]
\newtheorem{df}[theorem]{\bf Definition}
\newtheorem{thm}[theorem]{\bf Theorem}
\newtheorem{cor}[theorem]{\bf Corollary}
\newtheorem{lem}[theorem]{\bf Lemma}

\newtheorem{prop}[theorem]{\bf Proposition}
\newtheorem{hypothesis}[theorem]{\bf Hypothesis}

\newtheorem{pf}{\it Proof:}

\newcommand{\qed}{\hfill\hbox{\rule{6pt}{6pt}}}

\makeatletter
\@addtoreset{equation}{section}
\makeatother

\title{\sc Existence of a ground state for the Nelson model with a singular perturbation}
\date{}
\author{Takeru Hidaka\\
Graduate School of Mathematics, Kyushu University\\
Fukuoka, Japan, 819-0395}
\begin{document}
\maketitle
\begin{abstract}
The existence of a ground state of the Nelson Hamiltonian with perturbations of the form $\sum_{j=1}^{4}c_{j}\phi^{j}$ with $c_{4}>0$ is considered.
The self-adjointness of the Hamiltonian and the existence of a ground state are proven for arbitrary values of coupling constants.
\end{abstract}
\section{Introduction}
The Nelson model introduced in \cite{N} describes $N$-quantum mechanical particles coupled to a scalar bose field.
Let $\omega$ be a boson dispersion relation which describes the energy of a single boson.
Then the free field Hamiltonian $H_{\mathrm{f}}$ is given by the second quantization of $\omega$:
\begin{eqnarray}
H_{\mathrm{f}}=d\Gamma(\omega).
\end{eqnarray}
Let $K=-\Delta+V$ be a Hamiltonian of a quantum mechanical particle. Then
the standard Nelson Hamiltonian is formally given by
\begin{eqnarray}
H_{\text{Nelson}}=K+H_{\mathrm{f}}+\alpha\phi(\rho).\label{i.0}
\end{eqnarray}
Here
$\alpha$ is a coupling constant and $\phi(\rho)$ a field operator smeared by a test function $\rho$.

We consider the Nelson model with $\phi(\rho)$ replaced by
the singular perturbation:
\begin{eqnarray}
P(\phi(\rho))=\sum_{j=1}^{4}c_{j}\phi(\rho)^{j} \label{i.2}
\end{eqnarray}
with $c_{4}>0$.
Thus the total Hamiltonian under consideration is
\begin{eqnarray}
H=K+H_{\mathrm{f}}+P(\phi(\rho)) \label{i.1}
\end{eqnarray}
with the domain $D(K)\cap D(H_{\mathrm{f}}) \cap D(\phi(\rho)^{4})$.
We suppose that $K$ has a compact resolvent, and $V_{-}^{1/2}$ is relatively bounded with respect to $(-\Delta)^{1/2}$, where $V_{-}\geq0$ is the negative part of $V$.

We are concerned with the spectrum of $H$ in the non-perturbative way.
The bottom of the spectrum of a Hamiltonian is called a ground state energy, and a eigenvector associated with the ground state energy a ground state.
We see that the bottom of the spectrum of $K+H_{\rm f}$ is equal to
the edge of continuum.
Then it is not trivial to show the existence of ground state of $H$ even when perturbations are not singular.

The main result of this paper is  to show (1) and (2) below:
\begin{itemize}
\item[(1)]
$H$ is self-adjoint and bounded from below;
\item[(2)]
$H$  has a ground state for all $\rho$ under some conditions.
\end{itemize}

({\bf Related models})
We review  here several models concerned  so far,
but this is incomplete list.

{\bf [Nelson model]}
Bach-Fr\"ohlich-Sigal \cite{bfs} show the existence and uniqueness of the
ground state of some general scalar model for sufficiently weak couplings.
This model includes the standard Nelson model.
Spohn \cite{Sp} proves however the existence of the ground state for the Nelson model for {\it arbitrary} values of
coupling constants but if $K$ has purely discrete spectrum.
G\'erard \cite{Ge} also shows the similar result, but the method is different from  \cite{Sp}.
Hiroshima and Sasaki \cite{hs}
shows the enhanced binding of the many body Nelson model, i.e.,the existence of ground states is shown
for sufficiently large couplings but
the existence of ground state of decoupled Hamiltonian is not assumed.

The results mentioned above are proven under the so called infrared regularity conditions.
Then the next task  is to study the case of no infrared regularity conditions.
Arai, Hirokawa and Hiroshima
\cite{ahh}
show the absence of ground state of
some abstract quantum field models without infrared regularity conditions. L\H{o}rinczi, Minlos and Sphon \cite{LMS},
Derezi\'nski  and G\'erard \cite{dg}, and Hirokawa \cite{hi} prove that the Nelson Hamiltonian has no ground states if the infrared regularity condition is not assumed.
Arai \cite{ar} shows however that the Nelson model without infrared regularity condition also has a ground state  if a non-Fock representation is taken.
See also \cite{ghps1,ghps2} for the Nelson model on a pseudo Riemannian manifold.

{\bf [The Pauli-Fierz model]}
The Pauli-Fierz model is a quantum field model in nonrelativistic quantum electrodynamics.
Its interaction is given by minimal coupling, and then
the spectral analysis turns to be hard due to the derivative coupling.
Bach, Fr\"ohlich and Sigal \cite{BFS} prove the existence of ground state for sufficiently weak couplings
but the infrared regularity condition is not assumed.
This is large difference between the Nelson model and the Pauli-Fierz model.
Griesemer, Lieb and Loss \cite{GLL}, and Lieb and Loss \cite{LL} show the existence of a ground state of the Pauli-Fierz Hamiltonian for arbitrary values of coupling constants under no infrared regularity condition.
 In \cite{GLL, LL}, the binding condition is introduced to show the existence of a ground state.
We extend this to the Pauli-Fierz model with a variable mass \cite{hidaka}.  This method is also applied  to the Nelson model by Sasaki \cite{Sa}.

{\bf[Singular perturbations]}
The model under consideration in this paper
is of the similar form of
the $(\phi^{4})_{2}$-model in quantum field theory.
This model describes bosons with self-interaction in two dimensional space-time.
Glimm and Jaffe \cite{GJ68,GJ69} considered the
spectral properties of the $(\phi^{4})_{2}$-model.
In this model, the dispersion relation is supposed to be strictly positive and the Hamiltonian is defined on a boson Fock space.
Miyao and Sasaki \cite{MS} show the existence of the ground state for a generalized spin-boson model with $\phi^{2}$-perturbation, and it is {not} supposed that
the particle Hamiltonian has a compact resolvent.
Takaesu \cite{T} shows the existence of a ground state for a generalized spin-boson model with a singular perturbation
of the form (\ref{i.2})
but for sufficiently small coupling constants.

 ({\bf Strategy})
As far as we know, it is new to show the existence of the ground state of (\ref{i.1}) for {\it all} values of a coupling constant.
Here we show an outline of our proofs.

By making use of \cite{A}, we can prove the essential self-adjointness of $H$.
First we show that $\phi^{4}$ is relatively bounded with respect to $H$.
This relative boundedness leads to the self-adjointness of (\ref{i.1}).

Next we show the existence of a ground state of $H$
by means of \cite{DG,Ge} for all values of coupling constants:
We define the Hamiltonian $H_{\sigma}$, with the test function $\rho$ replaced by $\rho_\sigma=\rho 1_{\{\sigma\leq \omega(k)\}}$,
and we show the existence of a ground state of $H_{\sigma}$ for all $\sigma > 0$.
We see that as $\sigma\to 0$, a normalized ground state of $H_{\sigma}$ weakly converges to a non-zero vector, which is then a normalized ground state of $H$.
To show this it is sufficient to show the boson number bound and the boson derivative bound of a normalized ground state of $H_{\sigma}$.
These are done in Lemmas \ref{l.4.3} and \ref{l.4.5}.
To show the boson derivative bound, we suppose the infrared regularity condition:
\begin{equation}
\label{l1}
\omega^{-5/4}\sup_{x\in\mathbb{R}^{d}}|\rho(x,\cdot)|\in L^{2}(\mathbb{R}_{k}^{d}).
\end{equation}
This infrared regularity condition is stronger than the standard infrared regularity condition:
\begin{equation}\label{l2}
\omega^{-1}\sup_{x\in\mathbb{R}^{d}}|\rho(x,\cdot)|\in L^{2}(\mathbb{R}_{k}^{d}).
\end{equation}
The condition (\ref{l1}) is used to show the
convergence:
$$\|
((E_\sigma-H_\sigma-\omega(k))^{-1}-
(E-H-\omega(k))^{-1})\rho_\sigma(k)P'(\phi_\sigma)\Phi_\sigma\|
\to 0$$ in $L^2$ as $\sigma\to 0$ in Lemma \ref{l.4.4}, where $\Phi_\sigma$ is a ground state of $H_\sigma$.
In the case of the standard Nelson model, 
$P'(\phi_\sigma)=1$. Then condition (\ref{l2}) is enough to show
this convergence. In the singular case, we need however (\ref{l1}) to control the upper bound of
$\|
((E_\sigma-H_\sigma-\omega(k))^{-1}-
(E-H-\omega(k))^{-1})\rho_\sigma(k)P'(\phi_\sigma)\Phi_\sigma\|$.

This paper is organized as follows:
 Section 2 is devoted to defining the Nelson Hamiltonian with a singular perturbation.
 In Section 3 we show the self-adjointness of $H$.
 In Section 4 we show the existence of a ground state of $H$ but with an infrared cutoff.
 Finally in Section 5 we show the existence of a ground state of $H$.
\section{Definition of the Nelson model with $P(\phi)$ perturbation}
\subsection{Preliminaries}
Here we introduce fundamental facts on Fock spaces and second quantizations.
Let $\mathcal{X}$ be a Hilbert space over the complex field ${\Bbb C}$.
Then
\begin{eqnarray}
\mathcal{F}_{\mathrm{b}}(\mathcal{X})=\bigoplus_{n=0}^{\infty}\left[\otimes_{s}^{n}\mathcal{X}\right]=\left\{\{\Psi^{(n)}\}_{n=0}^{\infty}\Big| \Psi^{(n)}\in \otimes_{s}^{n}\mathcal{X}, n\geq 0,\sum_{n=0}^{\infty}\Vert \Psi^{(n)}\Vert^{2}<\infty \right\}
\end{eqnarray}
is called the boson Fock space over $\mathcal{X}$, where $\otimes_{s}^{n}\mathcal{X}$ denotes the symmetric tensor product of $\mathcal{X}$ and $\otimes_{s}^{0}\mathcal{X}=\mathbb{C}$.
Let $\Omega=\{1,0,\cdots\}\in \mathcal{F}_{\mathrm{b}}(\mathcal{X})$ be the Fock vacuum.
The number operator $N$ is defined by
\begin{eqnarray}
(N\Psi)^{(n)}=n\Psi^{(n)}
\end{eqnarray}
with the domain
\begin{eqnarray}
D(N)= \left\{\{\Psi^{(n)}\}_{n=0}^{\infty}\in\mathcal{F}_{\mathrm{b}}(\mathcal{X})\;\Big|\;\sum_{n=0}^{\infty}n^{2} \Vert\Psi^{(n)}\Vert^{2}<\infty\right\}.
\end{eqnarray}
The finite particle subspace of $\mathcal{F}_{\mathrm{b}}(\mathcal{X})$ is a dense subspace of $\mathcal {F}$, which is  given by
\begin{eqnarray}
\mathcal{F}_{\mathrm{b},0}(\mathcal{X})
=\left\{\{\Psi^{(n)}\}_{n=0}^{\infty}\in\mathcal{F}_{\mathrm{b}}(\mathcal{X})\,|\,\Psi^{(n)}=0 \text{ except for finitely many $n$}\right\}.
\end{eqnarray}
The creation operator $a^{\dagger}(f)$ smeared by $f\in \mathcal{X}$ is also given by
\begin{eqnarray}
(a^{\dagger}(f)\Psi)^{(n)}=\sqrt{n}S_{n}(f\otimes \Psi^{(n-1)}),\;n\geq 1,
\end{eqnarray}
and $(a^{\dagger}(f)\Psi)^{(0)}=0$
with the domain
\begin{eqnarray}
D(a^{\dagger}(f))=\left\{\Psi\in\mathcal{F}_{\mathrm{b}}\;\Big|\; \sum_{n=1}^{\infty}\Vert\sqrt{n}S_{n}(f\otimes \Psi^{(n-1)})\Vert^{2}<\infty \right\}.
\end{eqnarray}
Here $S_{n}$ is the symmetrization operator on $\otimes^{n} \mathcal{X}$.
The annihilation operator smeared by $f\in\mathcal{X}$ is given by the adjoint of $a^{\dagger}(f)$:
\begin{eqnarray}
a(f)=(a^{\dagger}(f))^{*}.
\end{eqnarray}
Note that $a(f)$ is antilinear in $f$, while $a^{\dagger}(f)$ is linear in $f$.
We see that $a(f)\lceil_{\otimes_{s}^{n}\mathcal{X}}$ is bounded from $\otimes_{s}^{n}\mathcal{X}$ to $\otimes_{s}^{n-1}\mathcal{X}$ and
$a^{\dagger}(f)\lceil_{\otimes_{s}^{n}\mathcal{X}}$ from $\otimes_{s}^{n}\mathcal{X}$ to $\otimes_{s}^{n+1}\mathcal{X}$.
$a(f)$ and $a^{\dagger}(f)$ satisfy canonical commutation relations:
\begin{eqnarray}
[a(f),a^{\dagger}(g)]=(f,g),\quad [a(f),a(g)]=[a^{\dagger}(f),a^{\dagger}(g)]=0.
\end{eqnarray}
Let $\mathcal{D}$ be a dense subspace of $\mathcal{X}$.
Then
\begin{eqnarray}
\mathcal{F}_{\mathrm{b,fin}}(\mathcal{D})=\mathcal{L}\{ \Omega, a^{\dagger}(f_{1})\cdots a^{\dagger}(f_{n})\Omega\;|\;n\in\mathbb{N}, f_{j}\in \mathcal{D}, j=1,\cdots n\}
\end{eqnarray}
is also dense in $\mathcal{F}_{\mathrm{b}}(\mathcal{X})$, where $\mathcal{L}\{\cdots \}$ denotes the linear hull of $\{\cdots\}$.
The Sigal field smeared by $f\in\mathcal{X}$ is given by
\begin{eqnarray}
\phi(f)=\frac{1}{\sqrt{2}}(\,a(f)+a^{\dagger}(f)\,).
\end{eqnarray}
Let $\mathcal{D}$ be a dense subspace and $f\in\mathcal{D}$.
Then $\phi(f)$ is essentially self-adjoint on $\mathcal{F}_{\mathrm{b,fin}}(\mathcal{D})$.
The Sigal field satisfies the following commutation relation:
\begin{eqnarray}
[\phi(f),\phi(g)]=i\Im (f,g).
\end{eqnarray}
When $\mathcal{X}=L^{2}(\mathbb{R}^{d})$, let
\begin{eqnarray}
(a(k)\Psi)^{(n)}(k_{1},\cdots,k_{n})=\sqrt{n+1}\Psi^{(n+1)}(k,k_{1},\cdots,k_{n}).
\end{eqnarray}
If $\Psi\in D(N^{1/2})$, then $a(k)\Psi\in \mathcal{F}_{\mathrm{b}}(L^{2}(\mathbb{R}^{d}))$ for almost every $k\in\mathbb{R}^{d}$.
For $\Psi\in D(a(f))$, $$(a(f)\Psi)^{(n)}(k_{1},\cdots,k_{n})=\int_{\mathbb{R}^{d}} \overline{f(k)}(a(k)\Psi)^{(n)}(k,k_{1},\cdots,k_{n})dk$$ holds.

Let $\mathcal{X}$ and $\mathcal{Y}$ be Hilbert spaces, and $T$ a densely defined closable operator from $\mathcal{X}$ to $\mathcal{Y}$.
Then $\Gamma(T)$ is defined by
\begin{eqnarray}
\Gamma(T)=\bigoplus_{n=0}^{\infty}\otimes^{n} T\lceil_{\otimes_{s}^{n}\mathcal{X}}
\end{eqnarray}
with $\otimes^{0} T =1$.
If $\mathcal{X}=\mathcal{Y}$, the second quantization of $T$ is defined by
\begin{eqnarray}
d\Gamma(T)=\bigoplus_{n=0}^{\infty}T^{(n)},
\end{eqnarray}
where $T^{(0)}=0$ and
\begin{eqnarray}
T^{(n)}=\overline{\sum_{j=1}^{n}1\otimes\cdots\otimes 1\otimes \stackrel{j\mathrm{th}}{\breve{T}}\otimes 1\otimes \cdots \otimes 1\lceil_{\otimes_{s}^{n}\mathcal{X}}},\quad n\geq 1.
\end{eqnarray}
Here $\overline{S}$ denotes the closure of an operator $S$.
The number operator $N$ can be written as
\begin{eqnarray}
N=d\Gamma(1).
\end{eqnarray}

We define the unitary operator $U_{\mathcal{X},\mathcal{Y}}$ from $\mathcal{F}_{\mathrm{b}}(\mathcal{X}\oplus \mathcal{Y})$ to $\mathcal{F}_{\mathrm{b}}(\mathcal{X})\otimes \mathcal{F}_{\mathrm{b}}(\mathcal{Y})$ by
\begin{eqnarray}
\lefteqn{U_{\mathcal{X},\mathcal{Y}}a^{\dagger}(f_{1}\oplus 0)\cdots a^{\dagger}(f_{n}\oplus 0)a^{\dagger}(0\oplus g_{1})\cdots a^{\dagger}(0\oplus g_{n})\Omega} \nonumber \\
&=&a^{\dagger}(f_{1})\cdots a^{\dagger}(f_{n})\Omega\otimes a^{\dagger}(g_{1})\cdots a^{\dagger}(g_{n})\Omega . \label{unitary}
\end{eqnarray}
Let $T$ be a densely defined closable operator from $\mathcal{X}$ to $\mathcal{X}\oplus \mathcal{X}$.
Then the operator $\check\Gamma(T): \mathcal F(\mathcal{X})\rightarrow \mathcal F(\mathcal{X})\otimes \mathcal F(\mathcal{X})$ is defined by
\begin{eqnarray}
\check\Gamma(T)=U_{\mathcal{X},\mathcal{X}}\Gamma(T).
\end{eqnarray}
\subsection{The Nelson Hamiltonian with $P(\phi)$ perturbation}
In this paper the number of quantum mechanical particles is supposed to be one, but the spatial dimension $d$.
Let $\mathcal{K}=L^{2}(\mathbb{R}_{x}^{d})$ and $\mathcal{F}_{\mathrm{b}}=\mathcal{F}_{\mathrm{b}}(L^{2}(\mathbb{R}_{k}^{d}))$.
The Hilbert space of state space is given by
\begin{eqnarray}
\mathcal{H}=\mathcal{K}\otimes\mathcal{F}_{\mathrm{b}},
\end{eqnarray}
where $\mathcal{K}$ describes the state space of a quantum mechanical particle, and $\mathcal{F}_{\mathrm{b}}$ that of bosons.
Let $\omega$ be a boson dispersion relation.
We suppose that $\omega$ is a densely defined, non-negative multiplication operator on $L^{2}(\mathbb{R}_{k}^{d})$.
Further conditions on $\omega$ are given later.
The free field Hamiltonian is given by $d\Gamma(\omega)$. Let $K$ be a Hamiltonian of the quantum mechanical particle.
Then the decoupled Hamiltonian is given by
\begin{eqnarray}
H_{0}=K\otimes 1+ 1\otimes d\Gamma(\omega)
\end{eqnarray}
with the domain $D(H_{0})=D(K\otimes 1)\cap D(1\otimes d\Gamma(\omega))$.
In what follows, we denote $T\otimes 1$ and $ 1\otimes S$ by $T$ and $S$, respectively, for simplicity unless confusion arises.
Let us now define a field operator $\phi$ in $\mathcal{H}$.
Let $\rho(x,k)$ be a test function such that $\rho(x,\cdot)\in L^{2}(\mathbb{R}_{k}^{d})$ for each $x\in\mathbb{R}^{d}$.
Then we set
\begin{eqnarray}
\phi(\rho(x,\cdot))=\frac{1}{\sqrt{2}}\left(a(\rho(x,\cdot)+ a^{\dagger}(\rho(x,\cdot)\right).
\end{eqnarray}
$\phi(\rho(x,\cdot))$ is essentially self-adjoint for each $x\in\mathbb{R}^{d}$ on
$$\mathcal{F}_{\mathrm{b,fin}}= \mathcal{L}\{\Omega, a^{\dagger}(h_{1})\cdots a^{\dagger}(h_{n})\Omega |\, n\in\mathbb{N},\;f,h_{i}\in C_{\mathrm{c}}(\mathbb{R}_{k}^{d}),i=1,\cdots n\}.$$
Then the field operator $\phi$ is defined by the constant fiber direct integral of $\overline{\phi(\rho(x,\cdot))}$:
\begin{eqnarray}
\phi=\phi(\rho)=\int^{\oplus}_{\mathbb{R}^{d}}\overline{\phi(\rho(x,\cdot))}dx.
\end{eqnarray}
Let
\begin{eqnarray}
P(x)=x^{4}+c_{3}x^{3}+c_{2}x^{2}+c_{1}x,
\end{eqnarray}
where $c_{j}$, $j=1,2,3,$ are arbitrary real numbers.
Then the Nelson Hamiltonian with $P(\phi)$ perturbation is given by
\begin{eqnarray}
H=H_{0}+ P(\phi)
\end{eqnarray}
with the domain
$D(H)=D(H_{0})\cap D(\phi^{4})$.
\subsection{Hypotheses and main theorems}
To show the self-adjointness of $H$ and the existence of a ground state of $H$, we introduce the following hypotheses.
\begin{hypothesis}[Hypotheses of $K$]\label{h.3.1}
\begin{enumerate}
\item $K$ is given by
\begin{eqnarray}
K=-\Delta+V
\end{eqnarray}
with the domain $D(K)=D(-\Delta)\cap D(V)$.
Here $V$ is a real-valued multiplication operator, which describes an external potential.
\item There exist constants $0<a<1$ and $b>0$ so that for all $\Psi\in D(V_{-}^{1/2})$, $\Psi\in D(|p|)$ and
\begin{eqnarray}\label{h.3.1.1}
\Vert V_{-}^{1/2} \Psi\Vert^{2} \leq a \Vert |p| \Psi \Vert^{2}+ b\Vert\Psi \Vert^{2}.
\end{eqnarray}
Here $V_{-}(x)=\max\{0,-V(x)\}$ and $p=-i\nabla_{x}$.
\item $K$ is a non-negative, self-adjoint operator and has a compact resolvent.
\end{enumerate}
\end{hypothesis}
\begin{hypothesis}[Hypotheses of $\omega$]\label{h.3.2}
\begin{enumerate}
\item $\omega \in C( \mathbb{R}_{k}^{d};[0,\infty))$;
\item $ \nabla \omega\in L^{\infty}(\mathbb{R}_{k}^{d})$;
\item $\omega(k)=0$ if and only if $k=0$.
\end{enumerate}
\end{hypothesis}
\begin{df}
Let $\mathcal{X}$ be a Hilbert space.
$f\in L^{\infty}(\mathbb{R}^{d};\mathcal{X})$ is said to be weakly differentiable if
there exists $g\in L^{\infty}(\mathbb{R}^{d};\mathcal{X})$ such that for all $\Psi\in\mathcal{X}$ and $\varphi\in C_{\mathrm{c}}^{\infty}(\mathbb{R}^{d})$,
\begin{eqnarray}
\int_{\mathbb{R}^{d}} (\partial_{j}\varphi)(x)(\Psi, f(x))_{\mathcal{X}}dx=-\int_{\mathbb{R}^{d}}\varphi(x)(\Psi,g(x))_{\mathcal{X}}dx.
\end{eqnarray}
In this case, we denote $g(x)$ by $\partial_{j}f(x)$.
\end{df}
\begin{hypothesis}[Hypotheses of $\rho$]\label{h.3.4}
$x\mapsto \rho(x,\cdot)$ is an element of $L^{\infty}(\mathbb{R}_{x}^{d};L^{2}(\mathbb{R}_{k}^{d}))$ and weakly twice differentiable. Moreover,
for each $k\in\mathbb{R}_{k}^{d}$, $\rho(k)=\rho(\cdot,k)$ is a bounded operator on $L^{2}(\mathbb{R}_{x}^{d})$ such that
\begin{eqnarray}
\omega^{-1/2}\Vert \rho(\cdot) \Vert,\; \omega\Vert\rho(\cdot) \Vert,\; \omega^{-1/2}\Vert\nabla_{x} \rho(\cdot)\Vert,\; \Vert \nabla_{x}\rho(\cdot) \Vert \in L^{2}(\mathbb{R}_{k}^{d}).
\end{eqnarray}
\end{hypothesis}
\begin{hypothesis}[Infrared regularity condition]\label{h.3.5}
It holds that
\begin{eqnarray}
\omega^{-5/4} \Vert \rho(\cdot) \Vert \in  L^{2}(\mathbb{R}_{k}^{d}).
\end{eqnarray}
\end{hypothesis}
We denote $\Vert \omega^{l} \Vert\rho(\cdot)\Vert \Vert_{L^{2}(\mathbb{R}^{d}_{k})}$ by $\Vert \omega^{l}\rho \Vert$ for $-5/4\leq l\leq 1$.
Let
\begin{eqnarray}
\mathcal{H}_{\mathrm{fin}}=
\mathcal{L}\{f\otimes\Omega, f\otimes a^{\dagger}(h_{1})\cdots a^{\dagger}(h_{n})\Omega |\, n\in\mathbb{N},f\in D(K),\;h_{i}\in C_{\mathrm{c}}^{\infty}(\mathbb{R}_{k}^{d}),\, 1\leq i\leq n \}. \hspace{-1cm}\nonumber \\
\end{eqnarray}
Now let us state the main theorems.
\begin{theorem}\label{main.0}
Suppose Hypotheses \ref{h.3.1} and \ref{h.3.4}.
Then $H$ is self-adjoint and essentially self-adjoint on $\mathcal{H}_{\mathrm{fin}}$.
\end{theorem}
\begin{theorem}\label{main}
Suppose Hypotheses  \ref{h.3.1}, \ref{h.3.2}, \ref{h.3.4} and \ref{h.3.5}.
Then $H$ has a ground state.
\end{theorem}
\section{Self-adjointness of $H$}
The following proposition on the essential self-adjointness is known.
\begin{prop} \label{prop of s.a.} \cite{A}
Let $\mathfrak{H}=\bigoplus_{n=0}^{\infty} \mathfrak{H}_{n}$ be the direct sum of Hilbert spaces $\mathfrak{H}_{n}$, $n=0,1,2,\cdots,$ and
$$
\hat{\mathfrak{H}}=\{\{\Psi^{(n)}\}_{n=0}^{\infty}\in \mathfrak{H}| \Psi^{(n)}=0 \text{ for all but finitely many $n$}\}.
$$
The number operator in $\mathfrak{H}$ is defined by
\begin{eqnarray}
(N_{\mathfrak{H}}\Psi)^{(n)}=n\Psi^{(n)},\quad n=0,1,2,\cdots
\end{eqnarray}
with the domain
\begin{eqnarray}
D(N_{\mathfrak{H}})=\left\{\{\Psi^{(n)}\}_{n=0}^{\infty}\in \mathfrak{H}\,\Big|\,\sum_{n=0}^{\infty}n^{2}\Vert\Psi^{(n)}\Vert^{2}<\infty\right\}.
\end{eqnarray}
Let $A_{n},$ $n=1,2,\cdots,$ be  self-adjoint operators in $\mathfrak{H}_{n}$ and $B$ a symmetric operator in $\mathfrak{H}$.
Put $A=\oplus_{n=0}^{\infty}A_{n}$.
Let $P_{n}$ be the projection from $\mathcal{X}$ to $\mathfrak{H}_{n}\subset \mathfrak{H}$: $P_{n}=1_{\{n\}}(N_{\mathfrak{H}})$.
Suppose that
\begin{enumerate}
\item $A+B$ is bounded from below;
\item $\hat{\mathfrak{H}}\subset D(B)$ and there exists a constant $n_{0}\geq 0$ such that
$(P_{m}\Psi,B P_{n}\Psi)=0$ whenever $|m-n|\geq n_{0}$;
\item there exist a constant $c$ and a linear operator $L$ in $\mathfrak{H}$ such that $\mathrm{Ran}(L\lceil_{D(L)\cap P_{n}\mathfrak{H}})\subset P_{n}\mathfrak{H}$ and
$$|(\Theta,B\Psi)|\leq c\Vert L\Theta\Vert \,\Vert (N_{\mathfrak{H}}+1)^{2}\Psi\Vert.$$
\end{enumerate}
Then
$A+B$ is essentially self-adjoint on $D(A)\cap \hat{\mathfrak{H}}$.
\end{prop}
\begin{lem}\label{l.3.5}
Suppose Hypotheses \ref{h.3.1} and \ref{h.3.4}. Then
$H$ is an essentially self-adjoint on $\mathcal{H}_{\textrm{fin}}$.
\end{lem}
\begin{pf}
By Proposition \ref{prop of s.a.}, $H$ is essentially self-adjoint on $D(H_{0})\cap \mathcal{H}_{0}$, where $\mathcal{H}_{0}=\mathcal{K}\otimes \mathcal{F}_{0}$.
Thus it suffices to show that $\mathcal{H}_{\mathrm{fin}}$ is a core for $\overline{H}$.
Let $\Psi\in D(H_{0})\cap \mathcal{H}_{0}$. Then there exists a number $n_{0}$ so that for all $n\geq n_{0}$, $\Psi^{(n)}=0$.
Since $\mathcal{H}_{\mathrm{fin}}$ is a core for $H_{0}$ by Proposition \ref{2.7} in Appendix and $P(\phi)\lceil_{\mathcal{K}\otimes(\otimes_{s}^{n_{0}}\mathcal{F})}$ is a bounded operator,
it is seen that there exists a sequence $\{\Psi_{j}\}_{j=1}^{\infty} \subset \mathcal{H}_{\mathrm{fin}}$ such that $\Psi_{j}\to\Psi$ and $H\Psi_{j}\to H\Psi$. Therefore the lemma follows.\qed
\end{pf}
\begin{lem}\label{l.3.6}
Suppose Hypothesis \ref{h.3.4}.
Let $\phi'_{j} = -i\phi(\partial_{x,j} \rho)$. Then $\phi^{n}\Psi\in D(|p|^{2})$ for $\Psi\in\mathcal{H}_{\mathrm{fin}}$ and $n\in\mathbb{N}$, and
$[\phi,p_{j}]=\phi'_{j}$ follows on $\mathcal{H}_{\mathrm{fin}}$.
\end{lem}
\begin{pf}
Let
$$
\Phi=f\otimes a^{\dagger}(f_{1})\cdots a^{\dagger}(f_{n})\Omega, \quad \Psi=g\otimes a^{\dagger}(g_{1})\cdots a^{\dagger}(g_{n-1})\Omega,
$$
where $f$ and $f_{k}\in C_{\mathrm{c}}^{\infty}(\mathbb{R}^{d})$, $k=1,\cdots, n$, and $g\in D(K)$ and $g_{k}\in C_{\mathrm{c}}^{\infty}(\mathbb{R}^{d})$, $k=1,\cdots, n-1$.
It can be computed as
\begin{eqnarray}
\lefteqn{(p_{j}\Phi, \phi\Psi)}\nonumber \\
\!\!&=&\!\!\!\!\!\frac{i}{\sqrt{2}}\int_{\mathbb{R}^{d}} (\partial_{j}\bar{f})(x)g(x)
(a(\rho_{x})a^{\dagger}(f_{1})\cdots a^{\dagger}(f_{n})\Omega,\, a^{\dagger}(g_{1})\cdots a^{\dagger}(g_{n-1})\Omega)dx \nonumber \\
\!\!&=&\!\!\!\!\!\frac{i}{\sqrt{2}}\sum_{l=1}^{n}\int_{\mathbb{R}^{d}} (\partial_{j}\bar{f})(x)g(x)(f_{l},\rho_{x})dx (a^{\dagger}(f_{1})\cdots \widehat{a^{\dagger}(f_{l})}\cdots a^{\dagger}(f_{n})\Omega,\, a^{\dagger}(g_{1})\cdots a^{\dagger}(g_{n-1})\Omega).\nonumber \\
\end{eqnarray}
Here the symbol $\widehat{ }$ denotes omission.
Since $\rho(x,\cdot)$ is weakly differentiable, we see that
\begin{eqnarray}
(p_{j}\Phi, \phi\Psi)=(\phi\Phi,p_{j}\Psi)+(\Phi,\phi_{j}' \Psi)=(\Phi, (\phi p_{j}+\phi_{j}')\Psi).\label{l.3.3.1}
\end{eqnarray}
Thus we obtain that $\phi\Psi\in D(p_{j})$ and
\begin{eqnarray}
p_{j}\phi\Psi=(\phi p_{j}+\phi_{j}')\Psi. \label{l.3.3.2}
\end{eqnarray}
By a similar computation, (\ref{l.3.3.2}) holds for all $\Psi\in \mathcal{H}_{\mathrm{fin}}$.
In a similar way, we can also see that $\phi^{n}\Psi\in D(|p|^{2})$.\qed
\end{pf}
\begin{thm} \label{t.3.7}
Suppose Hypotheses \ref{h.3.1} and \ref{h.3.4}. Then
there exists a constant $C$ such that for all $\Psi\in D(\overline{H})$,
\begin{eqnarray}
\Vert \phi^{n}\Psi \Vert \leq C\Vert (\overline{H}+1) \Psi\Vert,\quad n=1,2,3,4. \label{t.3.7.0}
\end{eqnarray}
\end{thm}
\begin{pf}
It is enough to show for the case of $n=4$.
Let $\Psi\in \mathcal{H}_{\text{fin}}$. It holds that
\begin{eqnarray}
\Vert \phi^{4}\Psi \Vert^{2}&=&\left\Vert \left(H-H_{0}-\sum_{k=1}^{3}c_{k}\phi^{k}\right)\Psi\right\Vert^{2}\nonumber \\
&=&\Vert H\Psi\Vert^{2}-2\Re(\phi^{4}\Psi, H_{0}\Psi)\nonumber \\
&&\quad-2\sum_{k=1}^{3}\Re(\phi^{4}\Psi,c_{k}\phi^{k}\Psi)-\left\Vert\left(H_{0}+\sum_{k=1}^{3}c_{k}\phi^{k}\right)\Psi\right\Vert^{2}\nonumber \\
&=& \Vert H\Psi\Vert^{2}-\left(\Psi,[\phi^{2},[\phi^{2},H_{0}]]\Psi\right) -2\sum_{k=1}^{3}c_{k}(\phi^{4}\Psi,\phi^{k}\Psi)\nonumber \\
&&\quad
-2\left\Vert H_{0}^{1/2}\phi^{2}\Psi\right\Vert^{2}-\left\Vert\left(H_{0}+\sum_{k=1}^{3}c_{k}\phi^{k}\right)\Psi\right\Vert^{2}.\label{t.3.7.11}
\end{eqnarray}
Take a sufficiently small $\epsilon>0$.
Since $\phi^{k}$, $k=1,2,3,$ are infinitesimally small with respect to $\phi^{4}$, there exists a constant $C_{1,\epsilon}>0$ such that
\begin{eqnarray}
-2\sum_{k=1}^{3}c_{k}(\phi^{4}\Psi,\phi^{k}\Psi)\leq \epsilon \Vert\phi^{4}\Psi\Vert^{2}+C_{1,\epsilon}\Vert\Psi\Vert^{2}.\label{t.3.7.12}
\end{eqnarray}
Thus by (\ref{t.3.7.11}) and (\ref{t.3.7.12}), we have
\begin{eqnarray}
\Vert \phi^{4}\Psi \Vert^{2}&\leq&\frac{1}{1-\epsilon}\Big( \Vert H\Psi\Vert^{2}+C_{1,\epsilon}\Vert\Psi\Vert^{2}+| (\Psi,[\phi^{2},[\phi^{2},H_{0}]]\Psi)|\Big).
\end{eqnarray}
Thus in order to prove (\ref{t.3.7.0}), it sufficies to show that for sufficiently small $0<\eta$, there exists a constant $C_{\eta}$ so that
\begin{eqnarray}
| (\Psi,[\phi^{2},[\phi^{2},H_{0}]]\Psi)|\leq \eta\Vert \phi^{4}\Psi\Vert^{2}+ C_{\eta}\Vert (H+1)\Psi \Vert^{2}. \label{t.3.7.1}
\end{eqnarray}
By Proposition \ref{p.2.1} (3) in Appendix, we have
\begin{eqnarray}
|(\Psi, [\phi^{2},[\phi^{2},d\Gamma(\omega)]]\Psi)| \leq 4\Vert\omega^{1/2}\rho\Vert^{2} \Vert \phi\Psi\Vert^{2}
\leq \epsilon \Vert \phi^{4}\Psi \Vert^{2}+ C_{2,\epsilon} \Vert \Psi\Vert^{2} . \label{t.3.7.2}
\end{eqnarray}
Let us estimate $|(\Psi, [\phi^{2},[\phi^{2},K]]\Psi)|$.
By Lemma \ref{l.3.6}, it is seen that
\begin{eqnarray}
|(\Psi, [\phi^{2},[\phi^{2},p_{j}^{2}]]\Psi )|&\leq& 2\Vert\,[\phi^{2},p_{j}]\Psi\Vert^{2}+|(\Psi,\{p_{j}[\phi^{2},[\phi^{2},p_{j}]]+[\phi^{2},[\phi^{2},p_{j}]]p_{j}\}\Psi)| \nonumber \\
&\leq& 2\Vert (\phi'_{j}\phi+\phi\phi'_{j})\Psi \Vert^{2}+8 |\Im (\rho,\partial_{x,j}\rho)||\Re(\phi^{2}\Psi,p_{j}\Psi)|.  \label{t.3.7.6}
\end{eqnarray}
Since $[\phi,\phi_{j}']=\Im (\rho,\partial_{x,j}\rho)$, by using the Schwarz inequality,
we see that
\begin{eqnarray}
|(\Psi, [\phi^{2},[\phi^{2},p_{j}^{2}]]\Psi )|\leq \epsilon \Vert\phi^{4}\Psi\Vert^{2}+ C_{3,\epsilon}(\Vert \phi'_{j}\phi \Psi\Vert^{2}+\Vert |p|\Psi\Vert^{2}+\Vert\Psi\Vert^{2}).\label{t.3.7.3}
\end{eqnarray}
Let us estimate $\Vert \phi'_{j}\phi \Psi\Vert^{2}$ in (\ref{t.3.7.3}). It holds that
\begin{eqnarray}
\Vert \phi'_{j}\phi \Psi\Vert^{2}&\leq& C_{4}(\phi\Psi,(d\Gamma(\omega)+1)\phi\Psi) \nonumber \\
&=&  C_{4}\{ (\phi^{2}\Psi,(d\Gamma(\omega)+1)\Psi)+(\phi\Psi,-i \phi(i\omega\rho)\Psi)\}. \label{t.3.7.7}
\end{eqnarray}
By the Schwarz inequality again, we have
\begin{eqnarray}
\Vert \phi_{j}'\phi\Psi\Vert^{2} \leq  \epsilon \Vert (d\Gamma(\omega)+1)\Psi\Vert^{2} +C_{4,\epsilon}\Vert \phi^{2} \Psi\Vert^{2}
+\frac{C_{4}}{2} (\Vert  \phi(i\omega\rho)\Psi\Vert^{2}+ \Vert \phi \Psi \Vert^{2}).\label{t.3.7.8}
\end{eqnarray}
Since $\phi(i\omega \rho)$ and $\phi^{k},$ $k=1,2,$ are infinitesimally small with respect to $d\Gamma(\omega)$ and $\phi^{4}$, respectively,
we see that
\begin{eqnarray}
\Vert \phi_{j}'\phi\Psi\Vert^{2}\leq \epsilon\Vert\phi^{4}\Psi\Vert^{2}+3\epsilon \Vert d\Gamma(\omega)\Psi\Vert^{2}+C_{5,\epsilon}\Vert \Psi\Vert^{2}.\label{t.3.7.13}
\end{eqnarray}
Since
\begin{eqnarray}
\Vert d\Gamma(\omega)\Psi \Vert^{2}&\leq& 2\Vert H\Psi\Vert^{2} +2\Vert P(\phi)\Psi\Vert^{2} \nonumber \\
&\leq& 2\Vert H\Psi\Vert^{2}+(2+\epsilon) \Vert \phi^{4}\Psi\Vert^{2}+ C_{6,\epsilon}\Vert\Psi\Vert^{2}, \label{t.3.7.14}
\end{eqnarray}
$\Vert \phi'_{j}\phi \Psi \Vert$ can be estimated by (\ref{t.3.7.13}) as
\begin{eqnarray}
\Vert \phi'_{j}\phi \Psi\Vert^{2} \leq \epsilon(3\epsilon+7) \Vert\phi^{4}\Psi\Vert^{2}+ 6\epsilon \Vert H\Psi\Vert^{2}+ (C_{5,\epsilon}+ 3\epsilon C_{6,\epsilon})\Vert\Psi\Vert^{2}. \label{t.3.7.9}
\end{eqnarray}
Next we estimate $\Vert|p|\Psi\Vert^{2}$ in (\ref{t.3.7.3}).
By (\ref{h.3.1.1}),
\begin{eqnarray}
\Vert |p|\Psi\Vert^{2} = (\Psi,-\Delta\Psi)\leq (\Psi, K\Psi)+(\Psi,V_{-}\Psi)\leq (\Psi, K\Psi)+a \Vert |p|\Psi\Vert^{2}+b\Vert\Psi\Vert^{2},
\end{eqnarray}
with $0<a<1$ and $0<b$.
Thus it holds that
\begin{eqnarray}
\Vert |p|\Psi\Vert^{2} \leq \frac{1}{1-a}(\Psi, (H-P(\phi))\Psi)+\frac{b}{1-a}\Vert\Psi\Vert^{2}.
\end{eqnarray}
Since $|\overline{H}|^{1/2}$ and $|P(\phi)|^{1/2}$ are infinitesimally small with respect to $\overline{H}$ and $\phi^{4}$, respectively, we have
\begin{eqnarray}
\Vert |p|\Psi\Vert^{2}\leq \epsilon (\Vert \phi^{4}\Psi\Vert^{2}+\Vert H\Psi\Vert^{2})+C_{7,\epsilon}\Vert\Psi\Vert^{2}. \label{t.3.7.5}
\end{eqnarray}
Therefore by (\ref{t.3.7.3}), (\ref{t.3.7.9}) and (\ref{t.3.7.5}), for sufficiently small $\epsilon'>0$, we have
\begin{eqnarray}
|(\Psi,[\phi^{2},[\phi^{2},K]\Psi)|\!=\left|\sum_{j=1}^{3}(\Psi,[\phi^{2},[\phi^{2},p_{j}^{2}]\Psi)\right|
\leq \epsilon'(\Vert\phi^{4}\Psi\Vert^{2}\!+\!\Vert H \Psi\Vert^{2})+C_{8,\epsilon'}\Vert\Psi\Vert^{2}. \label{t.3.7.10}
\end{eqnarray}
Therefore (\ref{t.3.7.1}) follows from (\ref{t.3.7.2}) and (\ref{t.3.7.10}).
Thus (\ref{t.3.7.0}) is proven for $\Psi\in \mathcal{H}_{\mathrm{fin}}$.
Since $\mathcal{H}_{\mathrm{fin}}$ is a core for $\overline{H}$,
(\ref{t.3.7.0}) also holds for all $\Psi\in D(\overline{H})$ by the closedness of $\phi^{4}$.\qed
\end{pf}

{\it Proof of Theorem \ref{main.0}:}
By Lemma \ref{l.3.5}, it suffices to show that
\begin{eqnarray}
 H=\overline{H}. \label{m.0}
\end{eqnarray}
Let $\Psi\in D(\overline{H})$.
Since $\mathcal{H}_{\mathrm{fin}}$ is a core for $\overline{H}$, there exists a sequence $\{\Psi_{j}\}_{j=1}^{\infty}$ so that $\Psi_{j}\in \mathcal{H}_{\mathrm{fin}}$ and
\begin{eqnarray}
\lim_{j\to\infty}(\Vert\Psi_{j}-\Psi\Vert+\Vert \overline{H}(\Psi_{j}-\Psi) \Vert)=0.
\end{eqnarray}
By the bound $\Vert\phi^{4}\Psi\Vert\leq C\Vert (\overline{H}+1)\Psi\Vert$ , $\{\phi^{4}\Psi_{j}\}_{j=1}^{\infty}$ is a Cauchy sequence. Since $\phi^{4}$ is closed, $\Psi\in D(\phi^{4})$ holds.
Since
\begin{eqnarray}
\Vert H_{0}(\Psi_{j}-\Psi_{k})\Vert\leq \Vert H(\Psi_{j}-\Psi_{k})\Vert + \Vert P(\phi)(\Psi_{j}-\Psi_{k}) \Vert,
\end{eqnarray}
$\{H_{0}\Psi_{j}\}_{j=1}^{\infty}$ is also a Cauchy sequence. Then $\Psi\in D(H_{0})$ by the closedness of $H_{0}$.
Thus
$D(\overline{H})\subset D(H_{0})\cap D(\phi^{4})=D(H)$.
Therefore
$(\ref{m.0})$ is obtained. \qed
\section{Existence of a ground state of $\tilde{H}_{\sigma}$ and $H_{\sigma}$}
\subsection{The Nelson Hamiltonian with an infrared cutoff $\sigma$}
The field operator with an infrared cutoff is given by
\begin{eqnarray}
\phi_{\sigma}=\phi(\rho_{\sigma}), \quad \sigma>0,
\end{eqnarray}
where
\begin{eqnarray}
\rho_{\sigma}=\rho 1_{\{k|\sigma\leq \omega(k)\}}.
\end{eqnarray}
We define $H_{\sigma}$  by
\begin{eqnarray}
H_{\sigma}=H_{0}+P(\phi_{\sigma})
\end{eqnarray}
with the domain $D(H_{\sigma})=D(H_{0})\cap D(\phi_{\sigma}^{4})$.
By Theorem \ref{main.0}, $H_{\sigma}$ is self-adjoint.
\begin{lem}\label{l.3.9}
Suppose Hypotheses \ref{h.3.1} and \ref{h.3.4}. Then
$H_{\sigma}$ converges to $H$ as $\sigma\to 0$ in the norm resolvent sense:
\begin{eqnarray}
\lim_{\sigma\to 0}\Vert (H_{\sigma}-z)^{-1}-(H-z)^{-1}\Vert=0
\end{eqnarray}
for all $z\in\mathbb{C}\setminus\mathbb{R}$.
\end{lem}
\begin{pf}
By the bound
\begin{eqnarray}
\Vert \phi_{\sigma}^{n}\Psi\Vert\leq C\Vert (\overline{H_{\sigma}}+1)\Psi\Vert,\quad n=1,2,3,4,
\end{eqnarray}
we see that
\begin{eqnarray}
\Vert (d\Gamma(\omega)+1)(H_{\sigma}-z)^{-1}\Vert <C,
\end{eqnarray}
with some constant $C$.
Take arbitrary vectors $\Theta\in \mathcal{H}$ and $ \Psi\in \mathcal{H} $. Then
\begin{eqnarray}
\!\!\lefteqn{|(\Theta, (H_{\sigma}-z)^{-1}-(H-z)^{-1}\Psi)|}\nonumber \\
\!\!\!\!\!\!&=&\!\!\!|((H_{\sigma}-\bar{z})^{-1}\Theta,P(\phi)(H-z)^{-1}\Psi)
-(P(\phi_{\sigma})(H_{\sigma}-\bar{z})^{-1}\Theta,(H-z)^{-1}\Psi)|.
\end{eqnarray}
Since $[\phi,\phi_{\sigma}]=0$ and $D(H)\cup D(H_{\sigma})\subset D(d\Gamma(\omega))\subset D(\phi_{\sigma}^{2})\cap D(\phi^{2})$, it follows that
\begin{eqnarray}
\lefteqn{|((H_{\sigma}-\bar{z})^{-1}\Theta,\phi^{4}(H-z)^{-1}\Psi)
-(\phi_{\sigma}^{4}(H_{\sigma}-\bar{z})^{-1}\Theta,(H-z)^{-1}\Psi)|}\nonumber \\
&\leq& |((\phi^{2}-\phi_{\sigma}^{2})  (H_{\sigma}-\bar{z})^{-1}\Theta, \phi^{2}(H-z)^{-1}\Psi)| \nonumber \\
&&\;+|(\phi_{\sigma}^{2}(H_{\sigma}-\bar{z})^{-1}\Theta, (\phi^{2}-\phi_{\sigma}^{2})(H-z)^{-1}\Psi)|\nonumber \\
&\leq&\Vert \phi(\rho-\rho_{\sigma}) \phi(\rho+\rho_{\sigma})(H_{\sigma}-\bar{z})^{-1}\Vert \,\Vert \phi^{2}(H-z)^{-1}\Vert \,\Vert\Theta\Vert\Vert\Psi\Vert \nonumber \\
&&+\Vert \phi_{\sigma}^{2}(H_{\sigma}-\bar{z})^{-1}\Vert \phi(\rho-\rho_{\sigma})\phi(\rho+\rho_{\sigma})(H-z)^{-1}\Vert \Vert\Theta\Vert\Vert\Psi\Vert \nonumber \\
&\leq& C' (\Vert\omega^{-1/2}( \rho -\rho_{\sigma}) \Vert+\Vert\omega( \rho -\rho_{\sigma}) \Vert) \Vert \Theta\Vert \Vert \Psi\Vert,
\end{eqnarray}
with some constant $C'$.
Similarly we see that
\begin{eqnarray}
\lefteqn{|((H_{\sigma}-\bar{z})^{-1}\Theta,P(\phi)(H-z)^{-1}\Psi)
-(P(\phi_{\sigma})(H_{\sigma}-\bar{z})^{-1}\Theta,(H-z)^{-1}\Psi)|}\nonumber \\
&&\hspace{2cm}\leq  C'' (\Vert\omega^{-1/2}( \rho -\rho_{\sigma}) \Vert+\Vert\omega (\rho -\rho_{\sigma})\Vert ) \Vert \Theta\Vert \Vert \Psi\Vert. \hspace{2cm}
\end{eqnarray}
with some constant $C''$.
Thus we obtain that
\begin{eqnarray}
\Vert (H_{\sigma}-z)^{-1}-(H-z)^{-1}\Vert\leq C'' (\Vert\omega^{-1/2}( \rho-\rho_{\sigma}) \Vert+\Vert\omega( \rho -\rho_{\sigma}) \Vert). \label{l.3.9.1}
\end{eqnarray}
Since the right hand side of (\ref{l.3.9.1}) converges to $0$ as $\sigma \to 0$,
the lemma follows.\qed
\end{pf}
We denote the ground state energies of $H_{\sigma}$ and $H$ by $E_{\sigma}$ and $E$, respectively: 
\begin{eqnarray}
E=\inf_{\Psi\in D(H),\Vert\Psi\Vert=1}(\Psi, H\Psi),\quad E_{\sigma}=\inf_{\Psi\in D(H_{\sigma}),\Vert\Psi\Vert=1}(\Psi, H_{\sigma}\Psi).
\end{eqnarray}
Since $H$, $H_{\sigma}\geq C$ with some constant $C$ independent of $\sigma$,
by Lemma \ref{l.3.9}, we obtain the following corollary:
\begin{cor}\label{c.3.13}
Suppose Hypotheses \ref{h.3.1} and \ref{h.3.4}.
Then
\begin{eqnarray}
\lim_{\sigma\to0}E_{\sigma}=E.
\end{eqnarray}
\end{cor}
Let us introduce a multiplication operator $\tilde{\omega}_{\sigma}$ below:
\begin{eqnarray}
&&\tilde{\omega}_{\sigma}\in C(\mathbb{R}^{d}),\quad \nabla\tilde{\omega}_{\sigma}\in L^{\infty}(\mathbb{R}_{k}^{d}),\\
&&\tilde{\omega}_{\sigma}(k)\geq \frac{\sigma}{2}\quad \text{for}\quad k\in\mathbb{R}^{d},\\
&&\tilde{\omega}_{\sigma}(k)=\omega(k)\quad \text{if}\quad |k|\geq \sigma.
\end{eqnarray}
Then we define the massive Hamiltonian $\tilde{H}_{\sigma}$ by
\begin{eqnarray}
\tilde{H}_{\sigma}=K\otimes 1+ 1\otimes d\Gamma(\tilde{\omega}_{\sigma})+ P(\phi_{\sigma}).
\end{eqnarray}
Similarly to the case of $H$ and $H_{\sigma}$, we can see that $\tilde{H}_{\sigma}$ is self-adjoint on $D(K)\cap D(d\Gamma(\tilde{\omega}_{\sigma}))\cap D(\phi_{\sigma}^{4})$.
\subsection{Extended Hamiltonian and existence of a ground state}
Throughout in this subsection, we suppose Hypotheses \ref{h.3.1}, \ref{h.3.2} and \ref{h.3.4}.
\begin{lem}\label{l.3.11}
\begin{enumerate}
\item Let $m\in\mathbb{Z}$. Then
$(N+1)^{-m}(\tilde{H}_{\sigma}-z)^{-1}(N+1)^{m+1}$
is a bounded operator and
\begin{eqnarray}
\Vert (N+1)^{-m}(\tilde{H}_{\sigma}-z)^{-1}(N+1)^{m+1}\Vert \leq C\sigma^{-1}(1+|\Im z|^{-1}) \label{l.3.1.1}
\end{eqnarray}
with some constant $C$ independent of $z$ and $\sigma$;
\item
Let $\chi \in C_{\mathrm{c}}^{\infty}(\mathbb{R})$. Then for all $l,m\in \mathbb{Z}$,
$(N+1)^{l}\chi (H)(N+1)^{m}$
is a bounded operator.
\end{enumerate}
\end{lem}
\begin{pf}
Let us show (1).
We denote $1_{\{n\}}(N)$ by $P_{n}$.
Let $m\geq 0$ and $\Psi\in D(N^{m+1})$.
Since $$P_{n}(H-z)^{-1}=\sum_{l=-4}^{4}P_{n}(H-z)^{-1}P_{n+l},$$ it follows that
\begin{eqnarray}
\lefteqn{\Vert (N+1)^{-m}(\tilde{H}_{\sigma}-z)^{-1}(N+1)^{m+1}\Psi \Vert^{2} } \nonumber \\
&=& \sum_{n=0}^{\infty} (n+1)^{-2m}\Vert P_{n}(\tilde{H}_{\sigma}-z)^{-1}(N+1)^{m+1} \Psi \Vert^{2} \nonumber \\
&\leq& C_{1}\sum_{n=0}^{\infty} (n+1)^{-2m}(n+5)^{2m}\sum_{l=-4}^{4}\Vert P_{n}(\tilde{H}_{\sigma}-z)^{-1}(N+1) P_{n+l} \Psi \Vert^{2} \nonumber \\
&\leq& C_{2} \Vert (\tilde{H}_{\sigma}-z)^{-1}(N+1) \Vert^{2}\sum_{n=0}^{\infty}\sum_{l=-4}^{4} \Vert P_{n+l} \Psi \Vert^{2} \nonumber \\
&\leq& C_{3}\sigma^{-2} (|\Im z|^{-1}+1)^{2} \Vert \Psi \Vert^{2},
\end{eqnarray}
where $C_{j}$, $j=1,2,3,$ are constants independent of $\sigma$ and $z$. In the last inequality, we used $\Vert N \Psi\Vert \leq \frac{2}{\sigma}\Vert d\Gamma(\tilde{\omega}_{\sigma})\Psi \Vert\leq \frac{C}{\sigma}\Vert (\tilde{H}_{\sigma}+1)\Psi \Vert$, since
$\tilde{\omega}_{\sigma}\geq \frac{\sigma}{2}$. Then (\ref{l.3.1.1}) follows.
When $m<0$, the lemma can be also proven similarly to the case of $m\geq 0$.
(2) can be proven similarly to \cite[Lemma 3.2 ii)]{DG}.  \qed
\end{pf}
Let us consider the extended Hilbert space defined by
\begin{eqnarray}
\mathcal{H}^{\text{ext}}=\mathcal{K}\otimes \mathcal{F}_{\mathrm{b}}\otimes \mathcal{F}_{\mathrm{b}}.
\end{eqnarray}
The decoupled Hamiltonian $\tilde{H}_{0,\sigma}$ is extended as
\begin{eqnarray}
\tilde{H}_{0,\sigma}^{\text{ext}}=\tilde{H}_{0,\sigma}\otimes 1_{\mathcal{F}_{\mathrm{b}}}+1_{\mathcal{H}}\otimes d\Gamma(\tilde{\omega}_{\sigma})
\end{eqnarray}
and
the total Hamiltonian $\tilde{H}_{\sigma}$ as
\begin{eqnarray}
\tilde{H}^{\text{ext}}_{\sigma}=\tilde{H}_{\sigma}\otimes 1_{\mathcal{F}_{\mathrm{b}}}+1_{\mathcal{H}}\otimes d\Gamma(\tilde{\omega}_{\sigma}).
\end{eqnarray}
Let us introduce a partition of unity such that
$j=(j_{0},j_{\infty})\in C^{\infty}(\mathbb{R}^{3};\mathbb{R}^{2})$,
$0\leq j_{0},j_{\infty}\leq 1$, $j_{0}^{2}+j_{\infty}^{2} = 1$ and
\begin{eqnarray}
j_{0}(x)=\left\{
\begin{array}{rl}
1 & \mbox{ if $|x| \leq 1$, }\\
0 & \mbox{ if $|x| \geq 2$.}
\end{array}
\right.
\end{eqnarray}
We set
$$j_{R}=(j_{0,R},j_{\infty,R})=(j_{0}(\cdot/R),j_{\infty}(\cdot/R))$$
and $$\hat{j}_{R}\Psi=(\hat{j}_{0,R}\Psi,\hat{j}_{\infty,R}\Psi)=(j_{0,R}(-i\nabla_{k})\Psi,\\j_{\infty,R}(-i\nabla_{k})\Psi).$$
Let us recall that $\check{\Gamma}(\hat{j}_{R}):\mathcal{F}_{\mathrm{b}}\rightarrow \mathcal{F}_{\mathrm{b}}\otimes \mathcal{F}_{\mathrm{b}}$ is defined by $U_{L^{2}(\mathbb{R}^{d}),L^{2}(\mathbb{R}^{d})}\Gamma(\hat{j}_{R})$.
\begin{lem}\label{l.3.13}
Let $\chi_{1}$, $\chi_{2}\in C_{\mathrm{c}}^{\infty}(\mathbb{R})$. Then
\begin{eqnarray}
\lim_{R\to 0}\left\Vert\left(\chi_{1}(\tilde{H}_{\sigma}^{\mathrm{ext}})\check{\Gamma}(\hat{j}_{R})-\check{\Gamma}(\hat{j}_{R})\chi_{1}(\tilde{H}_{\sigma})\right)\chi_{2}(\tilde{H}_{\sigma}) \right\Vert=0.
\end{eqnarray}
\end{lem}
\begin{pf}
By Helffer-Sj\"ostrand's formula, it is seen that
\begin{eqnarray}
\lefteqn{\left( \chi_{1}(\tilde{H}_{\sigma}^{\mathrm{ext}})\check{\Gamma}(\hat{j}_{R})-\check{\Gamma}(\hat{j}_{R})\chi(\tilde{H}_{\sigma}) \right)\chi_{2}(\tilde{H}_{\sigma})}\nonumber \\
&&=\frac{i}{2\pi}\int_{\mathbb{C}}\partial_{\bar{z}}\tilde{\chi}_{1}(z) (z-\tilde{H}_{\sigma}^{\mathrm{ext}})^{-1}(\tilde{H}_{\sigma}^{\mathrm{ext}}\check{\Gamma}(\hat{j}_{R})-\check{\Gamma}(\hat{j}_{R})\tilde{H}_{\sigma})(z-\tilde{H}_{\sigma})^{-1}\chi_{2}(\tilde{H}_{\sigma})dz d\bar{z}.\nonumber \\
\label{l.3.13.1}
\end{eqnarray}
Here $dzd\bar{z}=-2idxdy$, $\partial_{\bar{z}}=\frac{1}{2}(\partial_{x}+i\partial_{y})$ and
$\tilde{\chi}_{1}$ is an almost analytic extension of $\chi_{1}$, which satisfies
\begin{eqnarray}
\tilde{\chi}_{1}(x)&=&\chi_{1}(x),\quad x\in\mathbb{R},\\
\tilde{\chi}_{1}&\in& C_{\mathrm{c}}^{\infty}(\mathbb{C}),\\
|\partial_{\bar{z}}\tilde{\chi}_{1}(z)|&\leq& C_{n}|\Im z|^{n},\quad n\in\mathbb{N} .
\end{eqnarray}
Let us estimate the integrand in (\ref{l.3.13.1}).
$
\tilde{H}_{\sigma}^{\mathrm{ext}}\check{\Gamma}(\hat{j}_{R})-\check{\Gamma}(\hat{j}_{R})\tilde{H}_{\sigma}
$
is equal to
\begin{eqnarray}
\Big(\tilde{H}_{0,\sigma}^{\mathrm{ext}}\check{\Gamma}(\hat{j}_{R})-\check{\Gamma}(\hat{j}_{R})\tilde{H}_{0,\sigma}\Big)+\Big( (P(\phi_{\sigma})\otimes 1_{\mathcal{F}_{\mathrm{b}}})\check{\Gamma}(\hat{j}_{R})-\check{\Gamma}(\hat{j}_{R})P(\phi_{\sigma})\Big).\label{l.3.13.2}
\end{eqnarray}
The first term of (\ref{l.3.13.2}) can be estimated as
\begin{eqnarray}
\Vert (\tilde{H}_{0,\sigma}^{\text{ext}}\check{\Gamma}(\hat{j}_{R})-\check{\Gamma}(\hat{j}_{R})\tilde{H}_{0,\sigma}) (N+1)^{-1}\Vert
&=&\Vert  d{\Gamma}(\hat{j}_{R},((\tilde{\omega}_{\sigma}\oplus \tilde{\omega}_{\sigma})\hat{j}_{R}-\hat{j}_{R}\tilde{\omega}_{\sigma}) (N+1)^{-1}\Vert\nonumber \\
&\leq& \sqrt{\Vert [\tilde{\omega}_{\sigma},\hat{j}_{0,R}] \Vert^{2}+\Vert [\tilde{\omega}_{\sigma},\hat{j}_{\infty,R}] \Vert^{2}},
\end{eqnarray}
where $d\Gamma(\hat{j}_{R},\,(\tilde{\omega}_{\sigma}\oplus\tilde{\omega}_{\sigma})\hat{j}_{R}-\hat{j}_{R}\tilde{\omega}_{\sigma})$ is defined by
$$
(d\Gamma(\hat{j}_{R},(\tilde{\omega}_{\sigma}\oplus\tilde{\omega}_{\sigma})\hat{j}_{R}-\hat{j}_{R}\tilde{\omega}_{\sigma})\Psi)^{(n)}
=\sum_{l=1}^{n} \hat{j}_{R}\otimes\cdots\otimes \overbrace{((\tilde{\omega}_{\sigma}\oplus\tilde{\omega}_{\sigma})\hat{j}_{R}-\hat{j}_{R}\tilde{\omega}_{\sigma})}^{l\mathrm{th}}\otimes \cdots \otimes \hat{j}_{R}\Psi^{(n)}
$$
for $n\geq 1$ and
$$
(d\Gamma(\hat{j}_{R},(\tilde{\omega}_{\sigma}\oplus\tilde{\omega}_{\sigma})\hat{j}_{R}-\hat{j}_{R}\tilde{\omega}_{\sigma})\Psi)^{(0)}=0.
$$
Let us estimate commutators $[\tilde{\omega}_{\sigma},\hat{j}_{0,R}]$ and $[\tilde{\omega}_{\sigma},\hat{j}_{\infty,R}]$. Note that
\begin{eqnarray}
(f(-i\nabla)g)(k)=(2\pi)^{-d/2}\int_{\mathbb{R}^{d}}(\mathcal{F}f)(s)g(k+s)ds,
\end{eqnarray}
for $f\in C_{\mathrm{c}}^{\infty}(\mathbb{R}^{d})$ and $g\in C_{\mathrm{c}}^{\infty}(\mathbb{R}^{d})$. Here $\mathcal{F}f$ denotes the Fourier transformation of $f$.
Then
\begin{eqnarray}
\lefteqn{\Vert [\hat{j}_{0,R},\tilde{\omega}_{\sigma}]f\Vert_{L^{2}}^{2} }\nonumber \\
&=&(2\pi)^{-d}\int_{\mathbb{R}^{d}} \left|\int_{\mathbb{R}^{d}}(\mathcal{F}j_{0})(\xi)(\tilde{\omega}_{\sigma}(k+\xi/R)-\tilde{\omega}_{\sigma}(k))f(k+\xi/R)d\xi\right|^{2} dk \nonumber \\
&\leq& (2\pi)^{-d}  \Vert (\mathcal{F} j_{0}) \langle \cdot \rangle^{d+1}  \Vert_{L^{2}}^{2} (\Vert \nabla \tilde{\omega}_{\sigma}\Vert_{L^{\infty}} /R)^{2} \int_{\mathbb{R}^{d}} \int_{\mathbb{R}^{d}}
\langle \xi \rangle^{-2d}|f(k+\xi/R)|^{2}  dkd\xi \nonumber \\
&\leq& \frac{(2\pi)^{-d}}{R^{2}}  \Vert  (\mathcal{F} j_{0}) \langle \cdot \rangle^{d+1}  \Vert_{L^{2}}^{2} \Vert \nabla \tilde{\omega}_{\sigma}\Vert_{L^{\infty}}^{2} \Vert \langle \cdot \rangle^{-d} \Vert^{2}_{L^{2}}\,\Vert f\Vert_{L^{2}}^{2},
\end{eqnarray}
where $\langle \xi \rangle =\sqrt{1+\xi^{2}}$.
Thus
\begin{eqnarray}
\Vert [\hat{j}_{0,R},\tilde{\omega}_{\sigma}] \Vert =\frac{\mathrm{const.}}{R}.
\end{eqnarray}
Similarly,
\begin{eqnarray}
\Vert\, [\hat{j}_{\infty,R},\tilde{\omega}_{\sigma}] \,\Vert = \Vert [\hat{j}_{\infty,R}-1, \tilde{\omega}_{\sigma}] \Vert = \frac{\mathrm{const.}}{R},
\end{eqnarray}
since $j_{\infty,R}-1\in C^{\infty}_{\mathrm{c}}(\mathbb{R}^{d})$.
Thus it is seen that
\begin{eqnarray}
\lim_{R\to\infty}\Vert (\tilde{H}_{0,\sigma}^{\text{ext}}\check{\Gamma}(\hat{j}_{R})-\check{\Gamma}(\hat{j}_{R})\tilde{H}_{0,\sigma}) (N+1)^{-1}\Vert=0. \label{l.3.13.3}
\end{eqnarray}
Let us consider the second term of (\ref{l.3.13.2}). It can be computed as
\begin{eqnarray}
\lefteqn{\phi_{0,\sigma}^{4} \check{\Gamma}(\hat{j}_{R})- \check{\Gamma}(\hat{j}_{R}) \phi_{\sigma}^{4}} \nonumber \\
&=&\sum_{k=0}^{3}   \phi_{0,\sigma}^{3-k}[  (\phi_{0,\sigma} \check{\Gamma}(\hat{j}_{R})- \check{\Gamma}(\hat{j}_{R})\phi_{\sigma}  ]\phi_{\sigma}^{k} \nonumber \\
&=&\sum_{k=0}^{3}   \phi_{0,\sigma}^{3-k}\left[\left(\phi_{0}((1-\hat{j}_{0,R})\rho_{\sigma})-\phi_{\infty}(\hat{j}_{\infty,R}\rho_{\sigma}) \right)\check{\Gamma}(\hat{j}_{R}) \right]\phi_{\sigma}^{k}\label{l.3.13.4}
\end{eqnarray}
on $\mathcal{H}_{\mathrm{fin}}$. Here we write $\phi_{0}(f)$ and $\phi_{\infty}(f)$ for $\phi(f)\otimes 1_{\mathcal{F}_{\mathrm{b}}}$ and $1_{\mathcal{H}}\otimes \phi(f)$, respectively.
Note that
\begin{eqnarray}
\lim_{R\to\infty}\left\Vert(N_{0}+  N_{\infty})^{-3/2}\left[\left(\phi_{0}((1-\hat{j}_{0,R})\rho_{\sigma})-\phi_{\infty}(\hat{j}_{\infty,R}\rho_{\sigma}) \right)\check{\Gamma}(\hat{j}_{R}) \right](N+1)^{-2} \right\Vert=0,
\end{eqnarray}
where $N_{0}=N\otimes 1_{\mathcal{F}_{\mathrm{b}}}$ and $N_{\infty}=1_{\mathcal{H}}\otimes N$.
Then by (\ref{l.3.13.3}) and (\ref{l.3.13.4}),
\begin{eqnarray}
\lim_{R\to\infty}\Vert ( \tilde{H}_{\sigma}^{\mathrm{ext}}\check{\Gamma}(\hat{j}_{R})-\check{\Gamma}(\hat{j}_{R})\tilde{H}_{\sigma} )(N+1)^{-5/2} \Vert =0.\label{l.3.13.5}
\end{eqnarray}
By Lemma \ref{l.3.11}, the integrand of (\ref{l.3.13.1}) can be estimated as
\begin{eqnarray}
\lefteqn{
|\partial_{\bar{z}}\tilde{\chi}_{1}(z)|\Vert (z-\tilde{H}_{\sigma}^{\text{ext}})^{-1}(\tilde{H}_{\sigma}^{\text{ext}}\check{\Gamma}(\hat{j}_{R})-\check{\Gamma}(\hat{j}_{R})\tilde{H}_{\sigma}) (z-\tilde{H}_{\sigma})^{-1} \chi_{2}(\tilde{H}_{\sigma})\Vert } \nonumber \\
&\leq& |\partial_{\bar{z}}\tilde{\chi}_{1}(z)|\Vert (z-\tilde{H}_{\sigma}^{\mathrm{ext}})^{-1}\Vert \Vert( \tilde{H}_{\sigma}^{\mathrm{ext}}\check{\Gamma}(\hat{j}_{R})-\check{\Gamma}(\hat{j}_{R})\tilde{H}_{\sigma} )(N+1)^{-5/2} \Vert \times \nonumber \\
&&\quad \Vert (N+1)^{5/2}(z-\tilde{H}_{\sigma})^{-1}(N+1)^{-3/2}\Vert\,\Vert(N+1)^{3/2}\chi_{2}(\tilde{H}_{\sigma})\Vert \nonumber \\
&\leq& C \Vert (\tilde{H}_{\sigma}^{\mathrm{ext}}\check{\Gamma}(\hat{j}_{R})-\check{\Gamma}(\hat{j}_{R})\tilde{H}_{\sigma} )(N+1)^{-5/2}\Vert \left( 1+\frac{|\partial_{\bar{z}}\tilde{\chi}_{1}(z)|}{|\Im z |^{2}}\right),\label{l.3.13.6}
\end{eqnarray}
where $C$ is a constant independent of $z$ and $R$.
From (\ref{l.3.13.5}) and (\ref{l.3.13.6}), the lemma follows. \qed
\end{pf}
\begin{lem}\label{t.3.4}
Let $\tilde{E}_{\sigma}$ denote the ground state energy of $\tilde{H}_{\sigma}$.
Let $\chi\in C^{\infty}_{\mathrm{c}}(\mathbb{R})$. Suppose that $\mathrm{supp}\,\chi\subset (-\infty,\tilde{E}_{\sigma}+\sigma /2)$.
Then
$\chi(\tilde{H}_{\sigma})$ is a compact operator.
In particular, $\tilde{H}_{\sigma}$ has a ground state.
\end{lem}
\begin{pf}
First, let us show that $\Gamma(\hat{j}_{0,R}^{2})\chi(\tilde{H}_{\sigma})$ is a compact operator.
Since for each $n\in\mathbb{N}$,
$$\left\Vert\Gamma(\hat{j}_{0,R}^{2})\chi(\tilde{H}_{\sigma})- \sum_{k=0}^{n} 1_{\{k\}}(N)\Gamma(\hat{j}_{0,R}^{2})\chi(\tilde{H}_{\sigma})\right\Vert
\leq \frac{1}{n+1}\Vert\Gamma(\hat{j}_{0,R}^{2})N\chi(\tilde{H}_{\sigma})\Vert,$$
$ \sum_{k=0}^{n} 1_{\{k\}}(N)\Gamma(\hat{j}_{0,R}^{2})\chi(\tilde{H}_{\sigma})$ uniformly converges to $\Gamma(\hat{j}_{0,R}^{2})\chi(\tilde{H}_{\sigma})$ as $n$ goes to infinity. Then
it suffices to show that $1_{\{k\}}(N)\Gamma(\hat{j}_{0,R}^{2})\chi(\tilde{H}_{\sigma})$ is compact. Note that
$$
T_{1}=(K+1)^{-1/2}\otimes \Gamma(\hat{j}_{0,R}^{2})(d\Gamma(\tilde{\omega}_{\sigma})+1)^{-1/2}1_{\{k\}}(N)
$$ is compact and
$$
T_{2}=\left((K+1)^{1/2}\otimes(d\Gamma(\tilde{\omega}_{\sigma})+1)^{1/2}\right)\chi(\tilde{H}_{\sigma})
$$ is a bounded. Thus the claim is obtained since $1_{\{k\}}(N)\Gamma(\hat{j}_{0,R}^{2})\chi(\tilde{H}_{\sigma})=T_{1}T_{2}$.
Since $\mathrm{supp} \chi \subset (-\infty,\tilde{E}_{\sigma}+\sigma /2)$, we see that
\begin{eqnarray}
\chi (\tilde{H}_{\sigma}^{\text{ext}})=(1_{\mathcal{H}} \otimes P_{0})\chi(\tilde{H}_{\sigma}^{\mathrm{ext}}),\label{t.3.14.1}
\end{eqnarray}
where $P_{0}$ is the projection from $\mathcal{F}_{\mathrm{b}}$ to the subspace spanned by the Fock vacuum.
We also see that
\begin{eqnarray}
\check{\Gamma}(\hat{j}_{R})^{*}(1_{\mathcal{H}} \otimes P_{0}) \check{\Gamma}(\hat{j}_{R})=\Gamma(\hat{j}_{0,R}^{2}).\label{t.3.14.2}
\end{eqnarray}
We can suppose $\chi\geq 0$. Then by Lemma \ref{l.3.13} and (\ref{t.3.14.1}),
\begin{eqnarray}
\chi(\tilde{H}_{\sigma})&=&\check{\Gamma}(\hat{j}_{R})^{*}\check{\Gamma}(\hat{j}_{R})\chi^{1/2}(\tilde{H}_{\sigma})\chi^{1/2}(\tilde{H}_{\sigma}) \nonumber \\
&=&\check{\Gamma}(\hat{j}_{R})^{*} (1_{\mathcal{H}} \otimes P_{0})\chi^{1/2}(\tilde{H}_{\sigma}^{\text{ext}}) \check{\Gamma}(\hat{j}_{R})\chi^{1/2}(\tilde{H}_{\sigma})+o(R^{0}),
\end{eqnarray}
where $o(R^{0})$ is a bounded operator converging to $0$ as $R\to\infty$ in the uniform norm.
By Lemma \ref{l.3.13} again and (\ref{t.3.14.2}),
\begin{eqnarray}
\chi(\tilde{H}_{\sigma})&=&\check{\Gamma}(\hat{j}_{R})^{*}(1_{\mathcal{H}} \otimes P_{0}) \check{\Gamma}(\hat{j}_{R}) \chi(\tilde{H}_{\sigma})+o(R^{0}) \nonumber \\
&=&\check{\Gamma}(\hat{j}_{0,R}^{2})\chi(\tilde{H}_{\sigma})+o(R^{0}).
\end{eqnarray}
Since $\Gamma(\hat{j}_{0,R}^{2})\chi(\tilde{H}_{\sigma})$ is a compact operator, $\chi(\tilde{H}_{\sigma})$ is also compact. \qed
\end{pf}
\begin{lem}
$H_{\sigma}$ has a ground state.
\end{lem}
\begin{pf}
Let us consider the unitary operator (\ref{unitary}) with $\mathcal{X}=L^{2}(\{k|\omega(k)\geq \sigma\})$ and
$\mathcal{Y}=L^{2}(\{k|\omega(k)< \sigma\})$.
We denote $U_{\mathcal{X},\mathcal{Y}}$ by $U_{\sigma}$.
We see that
\begin{eqnarray}
U_{\sigma}\tilde{H}_{\sigma}U_{\sigma}^{*}=1\otimes d\Gamma(\tilde{\omega}_{\sigma})+H'_{\sigma}\otimes 1, \label{t.3.17.1}
\end{eqnarray}
where
\begin{eqnarray}
H'_{\sigma}=K+d\Gamma(\omega 1_{\{\omega(k)\geq \sigma\}})+P(\phi_{\sigma})
\end{eqnarray}
with the domain $D(K)\cap D(d\Gamma(\omega 1_{\{\omega(k)\geq \sigma\}}))\cap D(\phi_{\sigma}^{4})$.
$H'_{\sigma}$ is self-adjoint.
Since $\tilde{H}_{\sigma}$ has a ground state by Lemma \ref{t.3.4},
$H'_{\sigma}$ also has a ground state by Proposition \ref{2.7} in Appendix.
Since
\begin{eqnarray}
U_{\sigma}H_{\sigma}U_{\sigma}^{*}=1\otimes d\Gamma(\omega_{\sigma})+H'_{\sigma}\otimes 1,
\end{eqnarray}
and $H'_{\sigma}$ has a ground state,
$H_{\sigma}$ also has a ground state. \qed
\end{pf}
\section{Proof of the existence of a ground state}
Let $\Phi_{\sigma}$ be a normalized ground state of $H_{\sigma}$.
\begin{lem}[Pull-through formula]\label{l.4.2}
Suppose Hypotheses \ref{h.3.1}, \ref{h.3.2} and \ref{h.3.4}.
For almost every $ k\in\mathbb{R}^{d}$, we have
\begin{eqnarray}
a(k)\Phi_{\sigma}=\frac{1}{\sqrt{2}}(E_{\sigma}-H_{\sigma}-\omega(k))^{-1}\rho_{\sigma}(k)P'(\phi_{\sigma})\Phi_{\sigma}.
\end{eqnarray}
Here
\begin{eqnarray}
P'(x)=\frac{dP}{dx}(x).
\end{eqnarray}
\end{lem}
\begin{pf}
Since $\Phi_{\sigma}$ is a ground state of $H_{\sigma}$, for all $f\in C_{\mathrm{c}}^{\infty}(\mathbb{R}_{k}^{d})$ and $\Theta\in \mathcal{H}_{\mathrm{fin}}$,
\begin{eqnarray}
((H_{\sigma}-E_{\sigma})\Theta, a(f)\Phi_{\sigma})&=&([a^{\dagger}(f),H_{\sigma}-E_{\sigma}]\Theta, \Phi_{\sigma}) \nonumber \\
&=&\left(\left(-a^{\dagger}(\omega f)+\frac{1}{\sqrt{2}}(\rho_{\sigma},f)P'(\phi_{\sigma})\right)\Theta, \Phi_{\sigma}\right) \nonumber \\
&=&\left(\Theta,\left(-a(\omega f)+\frac{1}{\sqrt{2}}(f,\rho_{\sigma})P'(\phi_{\sigma})\right)\Phi_{\sigma}\right)\label{l.4.2.1}
\end{eqnarray}
follows. Since (\ref{l.4.2.1}) is equal to
\begin{eqnarray}
\int_{\mathbb{R}_{k}^{d}} \overline{f(k)}((E_{\sigma}-H_{\sigma}-\omega(k))\Theta, a(k)\Phi_{\sigma})dk=
\frac{1}{\sqrt{2}}\int_{\mathbb{R}_{k}^{d}}\overline{f(k)}(\Theta, \rho_{\sigma}(k)P'(\phi_{\sigma})\Phi_{\sigma})dk, \label{l.4.2.2}
\end{eqnarray}
it holds that for almost every $k\in\mathbb{R}^{d}$,
\begin{eqnarray}
((E_{\sigma}-H_{\sigma}-\omega(k))\Theta, a(k)\Phi_{\sigma})=\frac{1}{\sqrt{2}}(\Theta,  \rho_{\sigma}(k)P'(\phi_{\sigma})\Phi_{\sigma}). \label{l.4.2.3}
\end{eqnarray}
Thus
$a(k)\Phi_{\sigma}\in D(H_{\sigma})$ for almost every $k$ and
\begin{eqnarray}
(E_{\sigma}-H_{\sigma}-\omega(k)) a(k)\Phi_{\sigma}=\frac{1}{\sqrt{2}} \rho_{\sigma}(k)P'(\phi_{\sigma})\Phi_{\sigma}.
\end{eqnarray}
$E_{\sigma}-H_{\sigma}-\omega(k)<0$ for $k\neq 0$. Then $(E_{\sigma}-H_{\sigma}-\omega(k))^{-1}$ exists for $k\neq 0$.
Thus the lemma follows. \qed
\end{pf}
\begin{lem}\label{l.4.3}
Suppose Hypotheses \ref{h.3.1}, \ref{h.3.2}, \ref{h.3.4}, and
$\omega^{-1}\Vert\rho(\cdot)\Vert\in L^{2}(\mathbb{R}_{k}^{d})$.
Then
$\Phi_{\sigma}\in D(N^{1/2})$ and
\begin{eqnarray}
\sup_{0<\sigma\leq 1}\Vert N^{1/2}\Phi_{\sigma}\Vert <\infty.
\end{eqnarray}
\end{lem}
\begin{pf}
By Lemma \ref{l.4.2}, it follows that
\begin{eqnarray}
\Vert N^{1/2}\Phi_{\sigma}\Vert^{2}&=&\int_{\mathbb{R}_{k}^{d}}\Vert a(k)\Phi_{\sigma}\Vert^{2}dk\nonumber\\
&=&\frac{1}{2}\int_{\mathbb{R}_{k}^{d}}\Vert (E_{\sigma}-H_{\sigma}-\omega(k))^{-1}\rho_{\sigma}(k)P'(\phi_{\sigma})\Phi_{\sigma}\Vert^{2} dk\nonumber \\
&\leq&\frac{1}{2}\Vert P'(\phi_{\sigma})\Phi_{\sigma}\Vert^{2}\int_{\mathbb{R}_{k}^{d}} \Vert (E_{\sigma}-H_{\sigma}-\omega(k))^{-1}\rho_{\sigma}(k)\Vert^{2}dk \nonumber \\
&\leq&C\left(\sup_{0<\sigma\leq 1}E_{\sigma}^{2}+1\right) \left\Vert (1+\omega^{-1})\rho\right\Vert^{2}<\infty
\end{eqnarray}
with a constant $C$.
Thus the lemma follows.\qed
\end{pf}
\begin{lem} \label{l.4.1}
Suppose Hypotheses \ref{h.3.1}, \ref{h.3.2}, \ref{h.3.4} and \ref{h.3.5}. Then
it holds that
\begin{eqnarray}
|E_{\sigma}-E|\in o(\sigma^{3/4}). \label{l.4.1.0}
\end{eqnarray}
Here $\sigma^{-3/4}o(\sigma^{3/4})$ converges to $0$ as $\sigma\to +0$.
\end{lem}
\begin{pf}
Let $0<\sigma<\sigma'<1$.
Take a sequence $\{\Phi_{\sigma}^{j} \}_{j=1}^{\infty}\subset \mathcal{H}_{\mathrm{fin}}$ such that
$$ \lim_{j\to\infty}(\Vert \Phi_{\sigma'}^{j} - \Phi_{\sigma'}\Vert +\Vert H_{\sigma'}(\Phi_{\sigma'}^{j}-\Phi_{\sigma'})\Vert)=0.$$
Since there exist constants $C$ and $C'>0$ so that for all $j$ and $k\in\mathbb{N}$,
\begin{eqnarray}
\Vert \phi_{\sigma}^{2}(\Phi_{\sigma'}^{j}-\Phi_{\sigma'}^{k}) \Vert&\leq& C\Vert (d\Gamma(\omega)+1)(\Phi_{\sigma'}^{j}-\Phi_{\sigma'}^{k})\Vert \nonumber \\
&\leq& C(\Vert(H_{\sigma'}+1)(\Phi_{\sigma'}^{j}-\Phi_{\sigma'}^{k})\Vert+\Vert P(\phi_{\sigma'})(\Phi_{\sigma'}^{j}-\Phi_{\sigma'}^{k})\Vert)\nonumber \\
&\leq& C'\Vert (H_{\sigma'}+1)(\Phi_{\sigma'}^{j}-\Phi_{\sigma'}^{k})\Vert,
\end{eqnarray}
it is seen that
\begin{eqnarray}
\lim_{j\to\infty}\Vert\phi_{\sigma}^{2}(\Phi_{\sigma'}^{j}-\Phi_{\sigma'})\Vert=0\label{l.4.1.1}
\end{eqnarray}
by the closedness of $\phi_{\sigma}^{2}$.
Note that $\sup_{0<\sigma<1} |E_{\sigma}|<\infty $ by Corollary \ref{c.3.13}.
Then by (\ref{l.4.1.1}), it holds that
\begin{eqnarray}
 E_{\sigma}-E_{\sigma'}\!\!\!\!\!\!
&&\!\leq \liminf_{j\to\infty}\frac{(\Phi_{\sigma'}^{j},H_{\sigma}\Phi_{\sigma'}^{j})-(\Phi_{\sigma'}^{j},H_{\sigma'}\Phi_{\sigma'}^{j})}{\Vert \Phi_{\sigma'}^{j}\Vert^{2}} \nonumber \\
&&\!= ( \phi_{\sigma}^{2}\Phi_{\sigma'}, (\phi_{\sigma}^{2}-\phi_{\sigma'}^{2}) \Phi_{\sigma'})+(\phi_{\sigma'}^{2}\Phi_{\sigma'},(\phi_{\sigma}^{2}-\phi_{\sigma'}^{2}) \Phi_{\sigma'})\nonumber \\
&&\!\quad +c_{3}\{ (\phi_{\sigma}\Phi_{\sigma'}, (\phi_{\sigma}^{2}-\phi_{\sigma'}^{2}) \Phi_{\sigma'})+(\phi_{\sigma'}^{2}\Phi_{\sigma'},(\phi_{\sigma}-\phi_{\sigma'}) \Phi_{\sigma'})\}\nonumber \\
&&\!\quad +c_{2}\{ (\phi_{\sigma}\Phi_{\sigma'}, (\phi_{\sigma}-\phi_{\sigma'}) \Phi_{\sigma'})+(\phi_{\sigma'}\Phi_{\sigma'},(\phi_{\sigma}-\phi_{\sigma'}) \Phi_{\sigma'})\}\nonumber \\
&&\!\quad +c_{1} (\Phi_{\sigma'}, (\phi_{\sigma}-\phi_{\sigma'}) \Phi_{\sigma'})\nonumber
\end{eqnarray}
\begin{eqnarray}
&&\;\quad\leq C_{1} \left(\left\Vert\frac{\rho_{\sigma'}-\rho_{\sigma}}{\sqrt{\omega}} \right\Vert+ \Vert \omega (\rho_{\sigma'}-\rho_{\sigma}) \Vert \right) \left(\left\Vert\frac{\rho_{\sigma}}{\sqrt{\omega}}\right\Vert+\Vert\omega\rho_{\sigma}\Vert \right)\left(\left\Vert\frac{\rho_{\sigma'}}{\sqrt{\omega}}\right\Vert+\Vert\omega\rho_{\sigma'}\Vert\right) \nonumber \\
&&\!\;\quad\quad \times \Vert (d\Gamma(\omega)+1)\Psi_{\sigma'} \Vert\nonumber \\
&&\;\quad\leq C_{2}\left( \left\Vert\frac{\rho_{\sigma'}-\rho_{\sigma}}{\sqrt{\omega}} \right\Vert+\Vert\omega (\rho_{\sigma'}-\rho_{\sigma})\Vert \right),\label{l.4.1.2}
\end{eqnarray}
with constants $C_{1}$ and $C_{2}$.
Note that
\begin{eqnarray}
\left|\frac{\rho_{\sigma'}(x,k)-\rho_{\sigma}(x,k)}{\sqrt{\omega (k)}}\right|\leq \sigma'^{3/4}\left|\frac{ \rho_{\sigma'}(x,k)-\rho_{\sigma}(x,k) }{ \omega^{5/4}(k) }\right| \label{a}
\end{eqnarray}
and
\begin{eqnarray}
|\omega(k) (\rho_{\sigma'}(x,k)-\rho_{\sigma}(x,k))|\leq \sigma'^{3/4} |\omega^{1/4}(k)(\rho_{\sigma'}(x,k)-\rho_{\sigma}(x,k) )|. \label{b}
\end{eqnarray}
Then by (\ref{l.4.1.2}), it is obtained that
\begin{eqnarray}
E_{\sigma}-E_{\sigma'}\leq  C\sigma'^{3/4}\left( \left\Vert\frac{\rho_{\sigma'}-\rho_{\sigma}}{\omega^{5/4}} \right\Vert+\Vert\omega^{1/4} (\rho_{\sigma'}-\rho_{\sigma})\Vert \right). \label{l.4.1.3}
\end{eqnarray}
Replacing $\sigma$ and $\sigma'$, we have
\begin{eqnarray}
|E_{\sigma'}-E_{\sigma}|\leq  C\sigma'^{3/4}\left( \left\Vert\frac{\rho_{\sigma'}-\rho_{\sigma}}{\omega^{5/4}} \right\Vert+\Vert\omega^{1/4} (\rho_{\sigma'}-\rho_{\sigma})\Vert \right).\label{l.4.1.4}
\end{eqnarray}
Taking $\sigma'\to 0$ on both sides of (\ref{l.4.1.4}), we obtain (\ref{l.4.1.0}),
since $\omega^{-5/4}\Vert\rho(\cdot)\Vert$, $\omega^{1/4}\Vert\rho(\cdot)\Vert\in L^{2}(\mathbb{R}^{d}_{k})$.
\end{pf}
\begin{lem} \label{l.4.4}
Suppose Hypotheses \ref{h.3.1}, \ref{h.3.2}, \ref{h.3.4} and \ref{h.3.5}. Then
\begin{eqnarray}
\lim_{\sigma\to 0}\int_{\mathbb{R}_{k}^{d}}\left\Vert a(k)\Phi_{\sigma}-\frac{1}{\sqrt{2}}(E-H-\omega(k))^{-1}\rho(k)P'(\phi_{\sigma})\Phi_{\sigma}\right\Vert^{2}dk=0.
\end{eqnarray}
\end{lem}
\begin{pf}
Applying the pull through formula, Lemma \ref{l.4.2}, we have
\begin{eqnarray}
\lefteqn{
a(k)\Phi_{\sigma}-\frac{1}{\sqrt{2}}(E-H-\omega(k))^{-1}P'(\phi_{\sigma})\rho_{\sigma}(k)\Phi_{\sigma}
}\nonumber \\
&=&-\frac{1}{\sqrt{2}}(E-H-\omega(k))^{-1}(\rho(k)-\rho_{\sigma}(k))P'(\phi_{\sigma})\Phi_{\sigma}\nonumber \\
&&\quad+\frac{1}{\sqrt{2}}\left((E_{\sigma}-H_{\sigma}-\omega(k))^{-1}-(E-H-\omega(k))^{-1}\right)\rho_{\sigma}(k)P'(\phi_{\sigma})\Phi_{\sigma}. \label{l.4.4.1}
\end{eqnarray}
First let us consider the first term of the right hand side of (\ref{l.4.4.1}).
By the bound $\Vert \phi^{n}\Psi\Vert\leq C \Vert (H+1)\Psi \Vert$, $n=1,2,3,4,$ we see that
\begin{eqnarray}
\lefteqn{\int_{\mathbb{R}_{k}^{d}} \Vert(E-H-\omega(k))^{-1}(\rho(k)-\rho_{\sigma}(k))P'(\phi_{\sigma})\Phi_{\sigma}\Vert^{2} dk} \nonumber \\
&&\leq \left\Vert\omega^{-1}(\rho-\rho_{\sigma})\right\Vert^{2} \Vert P'(\phi_{\sigma})\Phi_{\sigma}\Vert^{2}
\leq C\left(\sup_{ 0< \sigma\leq 1} E_{\sigma}^{2}+1\right) \left\Vert\omega^{-1}(\rho-\rho_{\sigma})\right\Vert^{2}\label{l.4.4.2}
\end{eqnarray}
with some constant $C$.
Since $\omega^{-1}\Vert\rho(\cdot)\Vert\in L^{2}(\mathbb{R}^{d}_{k})$, we obtain that
\begin{eqnarray}
\lim_{\sigma\to0}\int_{\mathbb{R}^{d}} \Vert(E-H-\omega(k))^{-1}(\rho(k)-\rho_{\sigma}(k))P'(\phi_{\sigma})\Phi_{\sigma}\Vert^{2} dk =0.\label{l.4.4.3}
\end{eqnarray}
Next let us consider the second term of the right hand side of (\ref{l.4.4.1}). We see that
\begin{eqnarray}
\lefteqn{(\Theta, \left((E_{\sigma}-H_{\sigma}-\omega(k))^{-1}-(E-H-\omega(k))^{-1}\right)\rho_{\sigma}(k)P'(\phi_{\sigma})\Phi_{\sigma})} \nonumber \\
&=& (E-E_{\sigma})\left((E_{\sigma}-H_{\sigma}-\omega(k))^{-1}\Theta, (E-H-\omega(k))^{-1}\rho_{\sigma}(k)P'(\phi_{\sigma})\Phi_{\sigma}\right) \nonumber \\
&&+ \left(P(\phi_{\sigma})(E_{\sigma}-H_{\sigma}-\omega(k))^{-1}\Theta, (E-H-\omega(k))^{-1}\rho_{\sigma}(k)P'(\phi_{\sigma})\Phi_{\sigma}\right) \nonumber \\
&&- \left((E_{\sigma}-H_{\sigma}-\omega(k))^{-1}\Theta, P(\phi)(E-H-\omega(k))^{-1}\rho_{\sigma}(k)P'(\phi_{\sigma})\Phi_{\sigma}\right)
\end{eqnarray}
for all $\Theta\in\mathcal{H}$. It holds that
\begin{eqnarray}
\lefteqn{|\left(\phi_{\sigma}^{4}(E_{\sigma}-H_{\sigma}-\omega(k))^{-1}\Theta, (E-H-\omega(k))^{-1}\rho_{\sigma}(k)P'(\phi_{\sigma})\Phi_{\sigma}\right)} \nonumber \\
&&- \left((E_{\sigma}-H_{\sigma}-\omega(k))^{-1}\Theta, \phi^{4}(E-H-\omega(k))^{-1}\rho_{\sigma}(k)P'(\phi_{\sigma})\Phi_{\sigma}\right)| \nonumber \\
&=&|(\phi_{\sigma}^{2}(E_{\sigma}-H_{\sigma}-\omega(k))^{-1}\Theta, (\phi_{\sigma}^{2}-\phi^{2})(E-H-\omega(k))^{-1}\rho_{\sigma}(k)P'(\phi_{\sigma})\Phi_{\sigma}) \nonumber \\
&&+((\phi_{\sigma}^{2}-\phi^{2})(E_{\sigma}-H_{\sigma}-\omega(k))^{-1}\Theta,\phi^{2}(E-H-\omega(k))^{-1}\rho_{\sigma}(k)P'(\phi_{\sigma})\Phi_{\sigma})|\nonumber \\
&\leq& C(\Vert\omega^{-1/2}(\rho-\rho_{\sigma}) \Vert+ \Vert \omega ({\rho-\rho_{\sigma})}\Vert ) \left(1+\frac{1}{\omega(k)^{2}}\right) \Vert\rho_{\sigma}(k)\Vert \,\Vert \Theta \Vert, \label{l.4.4.10}
\end{eqnarray}
since
\begin{eqnarray}
\hspace{-2mm}\Vert(d\Gamma(\omega)+1)(E_{\sigma}-H_{\sigma}+1)^{-1}\Vert,\Vert(d\Gamma(\omega)+1)(E-H+1)^{-1}\Vert\leq C'\!\left(\!1+\frac{1}{\omega(k)^{2}}\!\right)\!\!,
\end{eqnarray}
where $C$ and $C'$ are constants.
Since $\Theta\in\mathcal{H}$ is arbitrary, in a similar way to (\ref{l.4.4.10}), we can see that
\begin{eqnarray}
\lefteqn{\Vert (\,(E_{\sigma}-H_{\sigma}-\omega(k))^{-1}-(E-H-\omega(k))^{-1})\rho_{\sigma}(k)P'(\phi_{\sigma})\Phi_{\sigma}\Vert} \nonumber \\
&\leq& C(\Vert\omega^{-1/2}(\rho-\rho_{\sigma}) \Vert+ \Vert \omega ({\rho-\rho_{\sigma})}\Vert + |E-E_{\sigma}| )\left(1+\frac{1}{\omega(k)^{2}}\right) \Vert\rho_{\sigma}(k)\Vert \nonumber \\
&\leq& C_{\sigma} \sigma^{3/4}  \left(1+\frac{1}{\omega(k)^{2}}\right) \Vert\rho_{\sigma}(k)\Vert \label{l.4.4.5}
\end{eqnarray}
for almost every $k\in\mathbb{R}^{d}$.
Here we used (\ref{a}) and (\ref{b}), and $C_{\sigma}$ is given by
\begin{eqnarray}
C_{\sigma}= \Vert \omega^{-5/4}(\rho-\rho_{\sigma}) \Vert+\Vert \omega^{1/4}(\rho-\rho_{\sigma})\Vert +\sigma^{-3/4}|E-E_{\sigma}|. \label{l.4.4.6}
\end{eqnarray}
Since $\sigma^{3/4}\Vert\rho_{\sigma}(k)\Vert\leq \omega(k)^{3/4}\Vert \rho_{\sigma}(k)\Vert$, by (\ref{l.4.4.5}) we see that
\begin{eqnarray}
\lefteqn{\Vert (\,(E_{\sigma}-H_{\sigma}-\omega(k))^{-1}-(E-H-\omega(k))^{-1})\,\rho_{\sigma}(k)P'(\phi_{\sigma})\Phi_{\sigma}\Vert}  \nonumber \\
&&\hspace{5cm}\leq C_{\sigma} \left(\omega(k)^{3/4}+\frac{1}{\omega(k)^{5/4}}\right) \Vert\rho_{\sigma}(k)\Vert.\label{l.4.4.7}
\end{eqnarray}
Since $|E-E_{\sigma}|=o(\sigma^{3/4})$, we have
\begin{eqnarray}
\lim_{\sigma\to 0}C_{\sigma}=0. \label{l.4.4.8}
\end{eqnarray}
Thus by (\ref{l.4.4.7}) and (\ref{l.4.4.8}), we obtain that
\begin{eqnarray}
\lim_{\sigma\to 0}\int_{\mathbb{R}_{k}^{d}}\Vert (\,(E_{\sigma}-H_{\sigma}-\omega(k))^{-1}-(E-H-\omega(k))^{-1})\rho_{\sigma}(k)P'(\phi_{\sigma})\Phi_{\sigma}\Vert^{2}dk
=0.\label{l.4.4.9}
\end{eqnarray}
Then we complete the lemma.
\end{pf}
\begin{lem}\label{l.4.5}
Suppose Hypotheses \ref{h.3.1}, \ref{h.3.2}, \ref{h.3.4} and \ref{h.3.5}.
Let $F\in C_{\mathrm{c}}^{\infty}(\mathbb{R}^{d})$ be such that $0\leq F(k)\leq 1$ for all $k\in\mathbb{R}^{d}$ and $F(0)=1$.
We set $F_{R}=F(\cdot/R)$ and $\hat{F}_{R}=F_{R}(D_{k})$, where $D_{k}=-i\nabla_{k}$.
Then
\begin{eqnarray}
\Vert d\Gamma(1-\hat{F}_{R})^{1/2}\Phi_{\sigma}\Vert=o(R^{0})+o(\sigma^{0}).
\end{eqnarray}
Here $o(R^{0})$ is a real number converging to $0$ as $R\to\infty$, and
$o(\sigma^{0})$ a real number converging to $0$ as $\sigma\to +0$.
\end{lem}
\begin{pf}
Since $\Phi_{\sigma}\in D(N^{1/2})$, $\Phi_{\sigma}\in D(d\Gamma(1-\hat{F}_{R})^{1/2})$. $F(D_{k})$ is defined as the bounded operator on $L^{2}(\mathbb{R}^{d};\mathcal{H})$ by
\begin{eqnarray}
(F(D_{k})\Psi(k))^{(n)}(x,k_{1},\cdots,k_{n})=F(D_{k})\Psi^{(n)}(k)(x,k_{1},\cdots,k_{n}),\quad 1\leq n,
\end{eqnarray}
and
$(F(D_{k})\Psi(k))^{(0)}(x,k_{1},\cdots,k_{n})=0$.
Then
\begin{eqnarray}
\Vert d\Gamma(1-\hat{F}_{R})^{1/2}\Phi_{\sigma}\Vert^{2}
= \int_{\mathbb{R}_{k}^{d}}(a(k)\Phi_{\sigma},(1-F(D_{k}/R))a(k)\Phi_{\sigma})dk. \label{l.4.5.1}
\end{eqnarray}
By the Schwarz inequality and Lemma \ref{l.4.4}, it is seen that
\begin{eqnarray}
(\ref{l.4.5.1})&\leq& \Vert N^{1/2}\Phi_{\sigma} \Vert \left( \int_{\mathbb{R}_{k}^{d}} \Vert (1-F(D_{k}/R))a(k)\Phi_{\sigma}\Vert^{2}dk\right)^{1/2}\nonumber \\
&=& \Vert N^{1/2}\Phi_{\sigma} \Vert \left(\int_{\mathbb{R}_{k}^{d}} \Vert (1-F(D_{k}/R))(E-H-\omega(k))^{-1}\rho(k)P'(\phi_{\sigma})\Phi_{\sigma}\Vert^{2}dk\right)^{1/2} \nonumber \\
&&\hspace{10.5cm}+o(\sigma^{0}). \label{l.4.5.2}
\end{eqnarray}
Let $\Theta\in L^{2}(\mathbb{R}_{k}^{d};\mathcal{H}_{\mathrm{fin}})$ with compact support and $T\in L^{2}(\mathbb{R}^{d}_{k};\mathcal{B}(\mathcal{H}))$. Then it holds that
\begin{eqnarray}
\lefteqn{\int(\Theta(k)^{(n)},F(D_{k})(T(k)\Psi)^{(n)})dk }\nonumber \\
&=&(2\pi)^{-d} \int\int \int\left( F(\xi)\Theta(s)^{(n)}e^{i\xi (k-s)},  (T(k)\Psi)^{(n)}\right)dsd\xi dk \nonumber \\
&=&(2\pi)^{-d/2} \int\int \left((\mathscr{F}F)(s-k) \Theta(s)^{(n)},(T(k)\Psi)^{(n)}\right)ds dk \nonumber \\
&=&(2\pi)^{-d/2} \int \left(\Theta(s)^{(n)}, (\mathbf{F}(D_{k})T)(s)\Psi^{(n)}\right) ds. \label{l.4.5.3}
\end{eqnarray}
Here $\mathscr{F}$ denotes the Fourier transformation and
$\mathbf{F}(D_{k})$ is the bounded operator on $L^{2}(\mathbb{R}^{d}_{k};\mathcal{B}(\mathcal{H}))$ defined by
\begin{eqnarray}
(\mathbf{F}(D_{k})T)(k)=(2\pi)^{-d/2}\int_{\mathbb{R}^{d}} (\mathscr{F}F)(s)T(k+s) ds. \label{l.4.5.4}
\end{eqnarray}
By (\ref{l.4.5.3}), we have
\begin{eqnarray}
F(D_{k})(T(k)\Psi)=(\mathbf{F}(D_{k})T)(k)\Psi. \label{l.4.5.5}
\end{eqnarray}
Then by (\ref{l.4.5.2}) and (\ref{l.4.5.5}), we see that
\begin{eqnarray}
\lefteqn{\Vert d\Gamma(1-\hat{F}_{R})^{1/2}\Phi_{\sigma}\Vert}\nonumber \\
&\leq& C \left(\int_{\mathbb{R}_{k}^{d}} \Vert (1-\mathbf{F}(D_{k}/R))(E-H-\omega(k))^{-1}\rho(k)\Vert^{2} dk\right)^{1/2}+o(\sigma^{0}).
\end{eqnarray}
Here $C$ is a constant independent of $\sigma$ and $R$. 
By Lemma \ref{l.4.7} below, the proof is complete. \qed
\end{pf}
\begin{lem}\label{l.4.7} \cite[Lemma 3.1]{Ge}
Suppose Hypotheses \ref{h.3.1}, \ref{h.3.2} and \ref{h.3.4} and suppose also that $\omega^{-1}\Vert \rho(\cdot) \Vert\in L^{2}(\mathbb{R}^{d}_{k})$.
Then it follows that
\begin{eqnarray}
\int_{\mathbb{R}^{d}} \Vert (1-\mathbf{F}(D_{k}/R))(E-H-\omega(k))^{-1}\rho(k)\Vert^{2} dk=o(R^{0}).
\end{eqnarray}
\end{lem}
{\it Proof of Theorem \ref{main}:}
Since $\Vert \Phi_{\sigma_{n}}\Vert =1$, we can take a sequence $\{\Phi_{\sigma_{n}}\}_{n=0}^{\infty}$ weakly converging to some vector $\Phi$ in $\mathcal{H}$:
\begin{eqnarray}
\text{w-}\lim_{n\to\infty}\Phi_{\sigma_{n}}= \Phi.
\end{eqnarray}
For all  $\Theta\in \mathcal{H}$ and $z\in \mathbb{C}\setminus \mathbb{R}$, it holds that
\begin{eqnarray}
(\Theta,(H_{\sigma_{n}}-z)^{-1}\Phi_{\sigma_{n}})=(\Theta,(E_{\sigma_{n}}-z)^{-1}\Phi_{\sigma_{n}}). \label{t.2.15.1}
\end{eqnarray}
Since $H_{\sigma_{n}}$ converges to $H$ in the norm resolvent sense, we see that by (\ref{t.2.15.1}) and Corollary \ref{c.3.13}
\begin{eqnarray}
(\Theta,(H-z)^{-1}\Phi)=(\Theta,(E-z)^{-1}\Phi ).
\end{eqnarray}
Since $\Theta$ is an arbitrary vector in $\mathcal{H}$, we have
\begin{eqnarray}
H\Phi=E\Phi. \label{t.2.15.2}
\end{eqnarray}
Thus $\Phi$ is a ground state of $H$ if and only if $\Phi\neq 0$.
We suppose $\Phi=0$.
Take $F\in C_{\mathrm{c}}^{\infty}(\mathbb{R}^{d})$ be such taht $0\leq F\leq 1$ and $F(0)=1$.
Since
$
\Gamma(\hat{F}_{R})1_{[0,\lambda]}(N)1_{[0,\lambda]}(H_{0})
$
is a compact operator, we see that
\begin{eqnarray}
\lim_{n\to\infty}\Vert\Gamma(\hat{F}_{R})1_{[0,\lambda]}(N)1_{[0, \lambda]}(H_{0})\Phi_{\sigma_{n}}\Vert=0. \label{t.2.15.3}
\end{eqnarray}
Note that
\begin{eqnarray}
\Vert (1- \Gamma(\hat{F}_{R}))\Psi \Vert\leq \Vert d\Gamma(1-\hat{F}_{R})^{1/2}\Psi\Vert \label{t.2.15.4}
\end{eqnarray}
 for all $\Psi\in D(d\Gamma(1-\hat{F}_{R})^{1/2})$.
Then we see that by (\ref{t.2.15.4}),
\begin{eqnarray}
\Vert \Phi_{\sigma_{n}}\Vert&\leq& \Vert \Gamma(\hat{F}_{R})\Phi_{\sigma_{n}} \Vert+ \Vert (1-\Gamma(\hat{F}_{R}))\Phi_{\sigma_{n}} \Vert \nonumber \\
&\leq& \Vert \Gamma(\hat{F}_{R})1_{[0,\lambda]}(N)1_{[0,\lambda]}(H_{0})\Phi_{\sigma_{n}} \Vert  +\Vert(1-1_{[0,\lambda]}(H_{0}))\Phi_{n} \Vert\nonumber \\
&&\quad +\Vert (1-1_{[0,\lambda]}(N))\Phi_{n} \Vert+ \Vert d\Gamma(1-\hat{F}_{R})^{1/2}\Phi_{\sigma_{n}}\Vert. \label{t.2.15.5}
\end{eqnarray}
Since $\sup_{n}(\Phi_{\sigma_{n}},N\Psi_{\sigma_{n}})<\infty$ and $\sup_{n}(\Psi_{\sigma_{n}},H_{0}\Phi_{\sigma_{n}})<\infty$, we have
\begin{eqnarray}
\Vert E_{N}((\lambda,\infty))\Phi_{n}\Vert,\; \Vert E_{H_{0}}((\lambda,\infty))\Phi_{n}\Vert \in O(\lambda^{-1}) . \label{t.2.15.6}
\end{eqnarray}
Here $E_{N}(\cdot)$ and $E_{H_{0}}(\cdot)$ are the spectral measures of $N$ and $H_{0}$, respectively.
Thus we see that
\begin{eqnarray}
\lim_{\lambda \to\infty}\sup_{n}\Vert (1-1_{[0,\lambda]}(N))\Phi_{n})\Vert=
\lim_{\lambda \to\infty}\sup_{n}\Vert (1-1_{[0,\lambda]}(H_{0}))\Phi_{n})\Vert=0. \label{t.2.15.7}
\end{eqnarray}
By (\ref{t.2.15.3}) and (\ref{t.2.15.7}), for an arbitrary $0<\epsilon<1$ we can take sufficiently large $0<\lambda$ so that
\begin{eqnarray}
\Vert \Phi_{\sigma_{n}}\Vert<\Vert 1_{[0,\lambda]}(N)1_{[0,\lambda]}(H_{0})\Gamma(\hat{F}_{R})\Phi_{\sigma_{n}} \Vert+\Vert d\Gamma(1-\hat{F}_{R})^{1/2}\Phi_{\sigma_{n}}\Vert
+\epsilon.\label{t.2.15.8}
\end{eqnarray}
Thus by (\ref{t.2.15.5}) and (\ref{t.2.15.8}),
\begin{eqnarray}
\limsup_{n\to\infty}\Vert \Phi_{\sigma_{n}}\Vert\leq \epsilon <1.
\end{eqnarray}
Since $\Phi_{\sigma_{n}}$ is a normalized vector in $\mathcal{H}$, this is a contradiction.
Therefore  $\Phi\neq 0$ and then
$\Phi$ is a ground state of $\mathcal{H}$. \qed
\section{Appendix}
Propositions \ref{p.2.1}-\ref{p.2.5} below are often used in this paper and well known. Let $\mathcal{X}$ and $\mathcal{Y}$ be Hilbert spaces.
\begin{prop}\label{p.2.1}\cite[Lemmas 2.7 and 2.8]{DG}
Let $T:\mathcal{X}\rightarrow \mathcal{Y}$ be a densely defined closable operator and $f\in D(T)$. Then
\begin{enumerate}
\item
\begin{eqnarray}
\Gamma(T)a^{\dagger}(f)=a^{\dagger}(Tf)\Gamma(T)
\end{eqnarray}
on $\mathcal{F}_{\mathrm{fin}}(D(T))$;
\item If $T$ is isometry, then
\begin{eqnarray}
\Gamma(T)a(f)=a(Tf)\Gamma(T)
\end{eqnarray}
on $\mathcal{F}_{\mathrm{fin}}(D(T))$;
\item If $\mathcal{X}=\mathcal{Y}$ and $f\in D(T)\cap D(T^{*})$, then
\begin{eqnarray}
[d\Gamma(T), a(f)]=-a(T^{*}f)\quad \text{and} \quad [d\Gamma(T), a^{\dagger}(f)]=a^{\dagger}(Tf)
\end{eqnarray}
on $\mathcal{F}_{\mathrm{b,fin}}(D(T))$.
\end{enumerate}
\end{prop}
\begin{prop}\cite[Proposition 8-6]{A00}
Let $\mathcal{X}=L^{2}(\mathbb{R}^{d})$.
\begin{enumerate}
\item Let $f$ be a function such that $0\leq f(k)<\infty$ for almost every $k$.
Then  $\Psi\in D(d\Gamma(f)^{1/2})$ if and only if
$$
\int_{\mathbb{R}^{d}}f(k)\Vert a(k)\Psi\Vert^{2}dk<\infty
$$
and in this case,
\begin{eqnarray}
\Vert d\Gamma(f)^{1/2}\Psi \Vert^{2}=\int_{\mathbb{R}^{d}}f(k)\Vert a(k)\Psi\Vert^{2}dk
\end{eqnarray}
holds.
Moreover if $f\in L^{2}(\mathbb{R}^{d})\cap L^{\infty}(\mathbb{R}^{d})$, it holds that
\begin{eqnarray}
\Vert d\Gamma(f(D))^{1/2}\Psi \Vert^{2}= \int_{\mathbb{R}^{d}} \sum_{n=1}^{\infty} \left((a(k)\Psi)^{(n)}, f(D_{k}) (a(k)\Psi)^{(n)}\right) dk
\end{eqnarray}
for all $\Psi\in D(d\Gamma(f(D))^{1/2})$. Here $D=-i\nabla$, and $D_{k}$ is the differential operator with respect to $k$.
\item
Let $f\in L^{2}(\mathbb{R}^{d})$, $\Phi\in\mathcal{F}_{\mathrm{b}}( L^{2}(\mathbb{R}^{d}))$ and $\Psi\in D(N^{1/2})$.
Then
\begin{eqnarray}
(\Phi,a(f)\Psi)=\int_{\mathbb{R}^{d}} \overline{f(k)}(\Phi, a(k)\Psi)dk.
\end{eqnarray}
\end{enumerate}
\end{prop}
\begin{prop}\cite[Proposition 4-24]{A00} and \cite[Lemma 2.1 i)]{DG}
\begin{enumerate}
\item Let $T$ be a self-adjoint operator with $\mathrm{ker}\, T=\{0\}$.
Suppose $f\in D(T^{-1/2})$.
Then for all $\Psi\in D(d\Gamma(T)^{1/2})$,
\begin{eqnarray}
\Vert a(f)\Psi \Vert&\leq& \Vert T^{-1/2}f \Vert \Vert d\Gamma(T)^{1/2}\Psi\Vert, \\
\Vert a^{\dagger}(f)\Psi \Vert^{2}&\leq& \Vert T^{-1/2}f\Vert^{2}  \Vert d\Gamma(T)^{1/2}\Psi\Vert^{2} +\Vert f\Vert^{2} \Vert\Psi \Vert^{2}.
\end{eqnarray}
\item
Let $l\in\mathbb{Z}$, $n\in\mathbb{N}$ and $f_{i}\in \mathcal{X}$, $i=1,\cdots, n$. Then
\begin{eqnarray}
\Vert (N+1)^{l}a^{\#}(f_{1})\cdots a^{\#}(f_{n}) (N+1)^{-l-\frac{n}{2}} \Vert \leq C_{n,l} \Pi_{i=1}^{n} \Vert f_{i}\Vert.
\end{eqnarray}
Here $a^{\#}(f)$ denotes $a(f)$ or $a^{\dagger}(f)$ and  $C_{n,l}$ is a constant depending on $n$ and $l$ but independent of $f_{i}$, $i=1,\cdots,n$.
\item
Let $T$ be a non-negative self-adjoint operator with $\mathrm{ker}\, T=\{0\}$.
Suppose that $f$, $g \in D(T)\cap D(T^{-1/2})$. Then
\begin{eqnarray}
\!\Vert a^{\#}(f)a^{\#}(g)\Psi \Vert \leq C \!\left(\left\Vert T^{-1/2}f\right\Vert+\Vert Tf\Vert\right)\left(\left\Vert T^{-1/2}g\right\Vert+\Vert T g\Vert\right)
\Vert (d\Gamma(T)+1)\Psi \Vert \nonumber \\
\end{eqnarray}
for $\Psi\in D(d\Gamma(T))$. Here $C$ is a constant independent of $T$, $f$, $g$ and $\Psi$.
\end{enumerate}
\end{prop}
\begin{prop}\label{2.7} \cite[Lemma 2-23, Corollary 2-27, Theorems 2-29 and 2-31]{A00}
Let $S$ and $T$ be non-negative self-adjoint operators in $\mathcal{X}$ and $\mathcal{Y}$ with cores
$\mathcal{D}_{1}$ and $\mathcal{D}_{2}$, respectively.
Then
\begin{enumerate}
\item $S\otimes 1_{\mathcal{Y}}$ and $1_{\mathcal{X}}\otimes T$ are strongly commuting;
\item
$S\otimes 1_{\mathcal{Y}}+ 1_{\mathcal{X}}\otimes  T$ is a self-adjoint operator and has a core $\mathcal{D}_{1}\hat{\otimes}\mathcal{D}_{2}$,
where $\hat{\otimes}$ denotes the algebraic tensor product;
\item It holds that for all $\Psi\in D(S\otimes 1_{\mathcal{Y}}+ 1_{\mathcal{X}}\otimes  T)$,
\begin{eqnarray}
\max\{\Vert (S\otimes 1_{\mathcal{Y}})\Psi\Vert, \Vert (1_{\mathcal{X}}\otimes T) \Psi\Vert\} \leq \Vert ( S\otimes 1_{\mathcal{Y}}+ 1_{\mathcal{X}}\otimes  T)\Psi \Vert;
\end{eqnarray}
\item For a densely defined closable operator $A$, we denote the spectrum of $A$ by $\sigma(A)$ and the point spectrum by $\sigma_{\mathrm{P}}(A)$, respectively. Then
\begin{eqnarray}
\sigma(S\otimes 1_{\mathcal{Y}}+ 1_{\mathcal{X}}\otimes  T)= \{\lambda+\mu|\lambda \in \sigma(S),\mu\in \sigma(T)\}
\end{eqnarray}
and
\begin{eqnarray}
\sigma_{\mathrm{P}}(S\otimes 1_{\mathcal{Y}}+ 1_{\mathcal{X}}\otimes  T)= \{\lambda+\mu|\lambda \in \sigma_{\mathrm{P}}(S),\mu\in \sigma_{\mathrm{P}}(T)\}.
\end{eqnarray}
\end{enumerate}
\end{prop}
\begin{prop}\cite[Theorem 4-55]{A00}\label{p.2.5}
\begin{enumerate}
\item
\begin{eqnarray}
U_{\mathcal{X},\mathcal{Y}}\mathcal{F}_{\mathrm{b,fin}}(\mathcal{X}\oplus\mathcal{Y})=\mathcal{F}_{\mathrm{b,fin}}(\mathcal{X})\hat{\otimes}\mathcal{F}_{\mathrm{b,fin}}(\mathcal{Y})
\end{eqnarray}
and
\begin{eqnarray}
U_{\mathcal{X},\mathcal{Y}}a^{\#}(f\oplus g)U_{\mathcal{X},\mathcal{Y}}^{-1}=a^{\#}(f)\otimes 1+1\otimes a^{\#}(g)
\end{eqnarray}
holds on $\mathcal{F}_{\mathrm{b,fin}}(\mathcal{X})\hat{\otimes}\mathcal{F}_{\mathrm{b,fin}}(\mathcal{Y})$.
\item \label{p.2.9}\cite[Theorem 4-56]{A00}
Let $T$ and $S$ be non-negative self-adjoint operators in $\mathcal{X}$ and $\mathcal{Y}$. Then
\begin{eqnarray}
U_{\mathcal{X},\mathcal{Y}}d\Gamma(T\oplus S)U_{\mathcal{X},\mathcal{Y}}^{-1}=d\Gamma(T)\otimes 1+1\otimes d\Gamma(S).
\end{eqnarray}
\end{enumerate}
\end{prop}
{\footnotesize
\section*{Acknowledgments}
I would like to thank Professor F. Hiroshima for his helpful comments and discussions.

\end{document}